\documentclass{amsart}
\usepackage{pb-diagram}
\usepackage{amssymb}
\usepackage[dvips]{epsfig}
\makeatletter
\def\arrowm#1#2{\def\dg@mx{#1}\def\dg@my{#2}}
\def\dgo@m{\let\dg@oldarrow=\dg@VECTOR\def\dg@VECTOR(##1,##2)##3{\put(\dg@mx,\dg@my){\dg@oldarrow(##1,##2){##3}}}}
\def\dg@delta{100}
\def\dgo@double{\let\dg@oldarrow=\dg@VECTOR\def\dg@VECTOR(##1,##2)##3{\put(0,\dg@delta){\dg@oldarrow(##1,##2){##3}}\put(0,-\dg@delta){\dg@oldarrow(##1,##2){##3}}}}
\makeatother

\newcommand\gfd{\mathcal D}
\newcommand\gfe{\mathcal E}
\newcommand\gff{\mathcal F}
\newcommand\gfg{\mathcal G}
\newcommand\gfh{\mathcal H}
\newcommand\gf{\mathcal X}
\newcommand\ball{\mathcal B}
\newcommand\free[1]{{\mathbb F_{#1}}}
\newcommand\C{{\mathbb C}}
\newcommand\R{{\mathbb R}}
\newcommand\Z{{\mathbb Z}} 
\newcommand\N{{\mathbb N}}

\newcommand\birth{\star}
\newcommand\death{\dagger}
\newcommand\bcount{\operatorname{bc}}
\newcommand\cbcount{\operatorname{cbc}}
\newcommand\tree{\mathcal T}
\newcommand\even{\operatorname{even}}
\newcommand\dist{\delta}
\newcommand\mF{\mathfrak F}
\newcommand\mG{\mathfrak G}
\newcommand\tens{\otimes}

\newcommand\Aut{\operatorname{Aut}}
\newcommand\idmatrix{\mathbb I}

\newtheorem{thm}{Theorem}[section]
\newtheorem{prop}[thm]{Proposition}
\newtheorem{lem}[thm]{Lemma}
\newtheorem{cor}[thm]{Corollary}
\theoremstyle{definition}\newtheorem{defn}[thm]{Definition}

\newcommand\emdef[1]{\emph{#1}
}
\newtheorem{problem}[thm]{Problem}

\begin{document}
\title{Counting Paths in Graphs}
\author{Laurent Bartholdi}
\date{\today}
\email{Laurent.Bartholdi@math.unige.ch}
\address{\parbox{.4\linewidth}{Section de Math\'ematiques\\ Universit\'e de Gen\`eve\\ CP 240, 1211 Gen\`eve 24\\ Switzerland}}
\nocite{coornaert-a:marches}
\begin{abstract}
  We give a simple combinatorial proof of a formula that extends a
  result by Grigorchuk~\cite{grigorchuk:phd,grigorchuk:rw} (rediscovered by
  Cohen~\cite{cohen:cogrowth}) relating cogrowth and spectral radius
  of random walks. Our main result is an explicit equation determining
  the number of `bumps' on paths in a graph: in a $d$-regular (not
  necessarily transitive) non-oriented graph let the series $G(t)$
  count all paths between two fixed points weighted by their length
  $t^{\text{length}}$, and $F(u,t)$ count the same paths, weighted as
  $u^{\text{number of bumps}}t^{\text{length}}$. Then one has
  $$\frac{F(1-u,t)}{1-u^2t^2} = \frac{G\left(\frac{t}{1+u(d-u)t^2}\right)}{1+u(d-u)t^2}.$$
  We then derive the circuit series of `free products' and `direct
  products' of graphs. We also obtain a generalized form of the
  Ihara-Selberg zeta function~\cite{bass:ihara,foata-z:bass}
\end{abstract}
\maketitle

\section{Introduction}
Let $\Gamma=\free{S}/N$ be a group generated by a finite set $S$,
where $\free{S}$ denotes the free group on $S$.  Let $f_n$ be the
number of elements of the normal subgroup $N$ of $\free{S}$ whose
minimal representation as words in $S\cup S^{-1}$ has length $n$; let
$g_n$ be the number of (not necessarily reduced) words of length $n$
in $S\cup S^{-1}$ that evaluate to $1$ in $\Gamma$; and let $d=|S\cup
S^{-1}|=2|S|$.  The numbers
$$\alpha=\limsup_{n\to\infty}\sqrt[n]{f_n},\qquad\nu=\frac1d\limsup_{n\to\infty}\sqrt[n]{g_n}$$
are called the \emph{cogrowth} and \emph{spectral radius} of
$(\Gamma,S)$.  The Grigorchuk Formula~\cite{grigorchuk:rw} states that
\begin{equation}\label{eq:grigformula}
  \nu = \begin{cases}
    \frac{\sqrt{d-1}}{d}\left(\frac\alpha{\sqrt{d-1}}+\frac{\sqrt{d-1}}\alpha\right) & \text{if }\alpha>\sqrt{d-1},\\
    \frac{2\sqrt{d-1}}{d} & \text{else}.
  \end{cases}
\end{equation}

We generalize this result to a somewhat more general setting: we
replace the group $\Gamma$ by a regular graph $\gf$, i.e.\ a graph
with the same number of edges at each vertex. Fix a vertex $\birth$ of
$\gf$; let $g_n$ be the number of circuits (closed sequences of edges)
of length $n$ at $\birth$ and let $f_n$ be the number of circuits of
length $n$ at $\birth$ with no backtracking (no edge followed twice
consecutively).  Then the same equation holds between the growth rates
of $f_n$ and $g_n$.

To a group $\Gamma$ with fixed generating set one associates its
Cayley graph $\gf$ (see Subsection~\ref{subs:gptheory}). $\gf$ is a
$d$-regular graph with distinguished vertex $\birth=1$; paths starting
at $\birth$ in $\gf$ are in one-to-one correspondence with words in
$S\cup S^{-1}$, and paths starting at $\birth$ with no backtracking
are in one-to-one correspondence with elements of $\free{S}$.  A
circuit at $\birth$ in $\gf$ is then precisely a word evaluating to
$1$ in $\Gamma$, and a circuit without backtracking represents
precisely one element of $N$.  In this sense results on graphs
generalize results on groups.  The converse would not be true: there
are even graphs with a vertex-transitive automorphism group that are
not the Cayley graph of a group~\cite{paschke:norm}.

Even more generally, we will show that, rather than counting circuits
and proper circuits (those without backtracking) at a fixed vertex, we
can count paths and proper paths between two fixed vertices and obtain
the same formula relating their growth rates.

These relations between growth rates are consequences of a stronger
result, expressed in terms of generating functions. Define the formal
power series
$$F(t)=\sum_{n=0}^\infty f_nt^n,\qquad G(t)=\sum_{n=0}^\infty g_nt^n.$$
Then assuming $\gf$ is $d$-regular we have
$$\frac{F(t)}{1-t^2} = \frac{G\left(\frac{t}{1+(d-1)t^2}\right)}{1+(d-1)t^2}.$$
This equation relates $F$ and $G$, and so relates \emph{a fortiori}
their radii of convergence, which are $1/\alpha$ and $1/(d\nu)$. We
re-obtain thus the Grigorchuk Formula.

Finally, rather than counting paths and proper paths between two fixed
vertices, we can count, for each $m\ge0$, the number of paths with $m$
backtrackings, i.e.\ with $m$ occurrences of an edge followed twice in
a row. Letting $f_{m,n}$ be the number of paths of length $n$ with $m$
backtrackings, consider the two-variable formal power series
$$F(u,t) = \sum_{m,n=0}^\infty f_{m,n}u^mt^n.$$
Note that $F(0,t)=F(t)$ and $F(1,t)=G(t)$. The following equation now
holds:
$$\frac{F(1-u,t)}{1-u^2t^2} = \frac{G\left(\frac{t}{1+u(d-u)t^2}\right)}{1+u(d-u)t^2}.$$
Letting $u=1$ in this equation reduces it to the previous one.

A generalization of the Grigorchuk Formula in a completely different
direction can be attempted: consider again a finitely generated group
$\Gamma$, and an exact sequence
$$1\longrightarrow\Xi\longrightarrow\Pi\longrightarrow\Gamma\longrightarrow1$$
where this time $\Pi$ is not necessarily free. Assume $\Pi$ is
generated as a monoid by a finite set $S$.  Let again $g_n$ be the
number of words of length $n$ in $S$ evaluating to $1$ in $\Gamma$,
and let $f_n$ be the number of elements of $\Xi$ whose minimal-length
representation as a word in $S$ has length $n$. Is there again a
relation between the $f_n$ and the $g_n$?  In
Section~\ref{sec:nonfree} we derive such a relation for $\Pi$ the
modular group $PSL_2(\Z)$.

Again there is a combinatorial counterpart; rather than considering
graphs one considers a locally finite cellular complex $\mathcal K$ such
that all vertices have isomorphic neighbourhoods. As before, $g_n$
counts the number of paths of length $n$ in the $1$-skeleton of
$\mathcal K$ between two fixed vertices; and $f_n$ counts elements of
the fundamental groupoid, i.e.\ homotopy classes of paths, between two
fixed vertices whose minimal-length representation as a path in the
$1$-skeleton of $\mathcal K$ has length $n$. We obtain a relation
between these numbers when $\mathcal K$ consists solely of triangles
and arcs, with no two triangles nor two arcs meeting; these are
precisely the complexes associated with quotients of the modular
group.

\section{Main Result}
Let $\gf$ be a graph, that may have multiple edges and loops. We make
the following typographical convention for the power series that will
appear: a series in the formal variable $t$ is written $G(t)$, or $G$
for short, and $G(x)$ refers to the series $G$ with $x$ substituted
for $t$.  Functions are written on the right, with $(x)f$ or $x^f$
denoting $f$ evaluated at $x$.

We start by the precise definition of graph we will use:
\begin{defn}[Graphs]
  A \emdef{graph} $\gf$ is a pair of sets $\gf=(V,E)$ and maps
  $$\alpha: E\to V,\qquad \omega: E\to V,\qquad \overline\cdot: E\to E$$
  satisfying
  $$\overline{\overline e} = e,\qquad \overline e^\alpha = e^\omega.$$
  The graph $\gf$ is said to be \emdef{finite} if both $E(\gf)$ and
  $V(\gf)$ are finite sets.

  A \emdef{graph morphism} $\phi:\gfg\to\gfh$ is a pair of maps
  $(V(\phi),E(\phi))$ with $V(\phi):V(\gfg)\to V(\gfh)$ and
  $E(\phi):E(\gfg)\to E(\gfh)$ and
  $$\overline{e E(\phi)} = \overline{e}E(\phi),\qquad e^\alpha V(\phi) = (eE(\phi))^\alpha.$$
  
  Given an edge $e\in E(\gf)$, we call $e^\alpha$ and $e^\omega$ $e$'s
  \emdef{source} and \emdef{destination}, respectively. We say two
  vertices $x,y$ are \emdef{adjacent}, and write $x\sim y$, if they
  are connected by an edge, i.e.\ if there exists an $e\in E(\gf)$
  with $e^\alpha=x$ and $e^\omega=y$. We say two edges $e,f$ are
  \emdef{consecutive} if $e^\omega=f^\alpha$. A \emdef{loop} is an
  edge $e$ with $e^\alpha=e^\omega$.

  The \emdef{degree} $\deg(x)$ of a vertex $x$ is the number
  of incident edges:
  $$\deg(x) = \#\{e\in E(\gf)|\, e^\alpha = x\} = \#\{e\in E(\gf)|\, e^\omega = x\}.$$
  If $\deg(x)$ is finite for all $x$, we say $\gf$ is \emdef{locally finite}.
  If $\deg(x)=d$ for all vertices $x$, we say $\gf$ is \emdef{$d$-regular}.
\end{defn}

Note that the involution $e\mapsto\overline e$ may have fixed points.
Even though the edges of $\gf$ are individually oriented, the graph
$\gf$ itself should be viewed as an non-oriented graph. In case
$\overline\cdot$ has no fixed point, $\gf$ can be viewed as a
geometric graph.

\begin{defn}[Paths]
  A \emdef{path} in $\gf$ is a sequence $\pi$,
  $$\pi=(v_0,e_1,v_1,e_2,\dots,e_n,v_n)$$
  of edges and vertices of $\gf$, with $e_i^\alpha=v_{i-1}$ and
  $e_i^\omega=v_i$ for all $i\in\{1,\dots,n\}$ and $n\ge0$. The
  \emdef{length} of the path $\pi$ is the number $n$ of edges in
  $\pi$. The \emdef{start} of the path $\pi$ is $\pi^\alpha = v_0$,
  and its \emdef{end} is $\pi^\omega = v_n$.  If $\pi^\alpha =
  \pi^\omega$, the path $\pi$ is called a \emdef{circuit} at
  $\pi^\alpha$. In most cases, we will omit the $v_n$ from the
  description of paths; they are necessary only if $|\pi|=0$, in which
  case a starting vertex must be specified. We extend the involution
  $\overline\cdot$ from edges to paths by setting
  $$\overline\pi = (v_n,\overline{e_n},\dots,v_1,\overline{e_1},v_0)$$
  (note that $\overline{\pi}$ is a path from $\pi^\omega$ to $\pi^\alpha$).

  We denote by $E^*(\gf)$ the set of paths, with a partially defined
  multiplication given by concatenation: if $\pi$ and $\rho$ be two
  paths with $\pi^\omega=\rho^\alpha$, their \emdef{product} is
  defined as $\pi\rho=(\pi_1,\dots,\pi_{|\pi|},\rho_1,\dots,\rho_{|\rho|})$.
  For two vertices $x,y\in V(\gf)$ we denote by $[x,y]$ the set of
  paths from $x$ to $y$. We turn $V(\gf)$ into a metric space by defining
  for vertices $x,y\in V(\gf)$ their \emdef{distance}
  $$\dist(x,y) = \min_{\pi\in[x,y]}|\pi|.$$
  The \emdef{ball} of radius $n$ at a vertex $x\in V(\gf)$ is the
  subgraph $\ball(x,n)$ of $\gf$ with vertex set
  $$V(\ball(x,n)) = \{y\in V(\gf):\,\dist(x,y)\le n\}$$
  and edge set
  $$E(\ball(x,n)) = \{e\in E(\gf):\,e^\alpha\in V(\ball(x,n))\}.$$
  We define $\alpha_\ball(e)=\alpha(e)$,
  $$\overline e_\ball=\begin{cases}\overline e&\text{ if }e^\omega\in
  E(\ball)\\ e\text{ else},\end{cases}$$
  and $\omega(e)=\alpha(\overline e)$.
\end{defn}
This definition amounts to ``wrapping around disconnected edges''. It
has the advantage of preserving the degrees of vertices.


\begin{defn}[Bumps, Labellings]
  We say a path $\pi$ has a \emdef{bump at $i$} if
  $\pi_i=\overline{\pi_{i+1}}$; if the location of the bump is
  unimportant we will just say $\pi$ has a bump. The \emdef{bump count}
  $\bcount(\pi)$ of a path $\pi$ is the number of bumps in $\pi$.
  A \emdef{proper path} in $\gf$ is a path $\pi$ with no bumps.
  
  Let $\Bbbk$ be a ring. A $\Bbbk$-\emdef{labelling} of the graph
  $\gf$ is a map
  $$\ell: E(\gf)\to\Bbbk.$$
  The simplest examples of labellings are:
  \begin{itemize}
  \item The \emdef{trivial labelling}, given by $\Bbbk=\Z$
    and~$e^\ell=1$ for all $e\in E(\gf)$;
  \item The \emdef{length labelling}, given by $\Bbbk=\Z[[t]]$
    and~$e^\ell=t$ for all $e\in E(\gf)$.
  \end{itemize}

  A $\Bbbk$-labelling $\ell$ of $\gf$ induces a map, still written
  $\ell:E^*(\gf)\to\Bbbk$, by setting
  $$\pi^\ell = \pi_1^\ell\pi_2^\ell\cdots\pi_n^\ell.$$
  
  The labelling $\ell:E^*(\gf)\to\Bbbk$ is \emdef{complete}~\cite[\S
  VI.2]{eilenberg:a} if for any vertex $x$ of $\gf$ and any set $A$ of
  paths in $\gf$ starting at $x$ there is an element $(A)\Sigma$ of
  $\Bbbk$, and this function $\Sigma$ satisfies
  $$(\{\pi\})\Sigma=\pi^\ell,\qquad(A\sqcup
  B)\Sigma=(A)\Sigma+(B)\Sigma$$ for all paths $\pi$ and disjoint sets
  $A$ and $B$ ($\sqcup$ denotes disjoint union).  If $A$ is infinite,
  it is customary, though abusive, to write $(A)\Sigma$ as
  $\sum_{\pi\in A}\pi^\ell$.
\end{defn}
If $\Bbbk$ is a topological ring ($\R$, $\C$, the formal power series
ring $\Z[[t]]$, etc.), completion of $\ell$ implies that
$\pi^\ell\to0$ when $|\pi|\to\infty$, but the converse does not hold.
The completeness condition becomes that
$$(A)\Sigma = \lim_{B\subset A, |B|<\infty}\sum_{\pi\in B}\pi^\ell$$
be a well-defined element of $\Bbbk$ for all $A$; i.e., the limit
exists. Generally, we define the following topology on $\Bbbk$: a
sequence $(A_i)\Sigma\in\Bbbk$ converges to $0$ if and only if
$\min_{\pi\in A_i}|\pi|$ tends to infinity.

In the sequel of this paper all labellings will be assumed to be
complete. The length labelling is complete for locally finite graphs;
more generally, $\ell$ is complete when $\Bbbk$ is a discretely valued
ring, $e^\ell$ has a positive valuation for all edges $e$, and $\gf$
is locally finite. An arbitrary ring $\Bbbk$ may be embedded in
$\Bbbk'=\Bbbk[[t]]$, where $t$ has valuation $1$ and $\Bbbk$ has
valuation $0$; if $\ell:E(\gf)\to\Bbbk$ is a labelling, we define
$\ell':E(\gf)\to\Bbbk'$ by $e^{\ell'}=te^\ell$; and $\ell'$ will be
complete as soon as $\gf$ is locally finite. In particular the length
labelling is obtained from the trivial labelling through this
construction.  In all the examples we consider the labelling is
defined in this manner.

Throughout the paper we will assume a graph $\gf$ and two vertices
$\birth,\death\in V(\gf)$ have been fixed. We wish to enumerate the
paths, counting their number of bumps, from $\birth$ to $\death$ in
$\gf$. For a given complete edge-labelling $\ell$, consider the series
$$\mG(\ell)=\sum_{\pi\in[\birth,\death]}\pi^\ell\in\Bbbk,\qquad
\mF(\ell)=\sum_{\pi\in[\birth,\death]}u^{\bcount(\pi)}\pi^\ell\in\Bbbk[[u]].$$
Beware that in general $\mG$ and $\mF$ also depend on the choice of
$\birth$ and $\death$.

For all vertices $x\in V(\gf)$ let
$$K_x = \left(1 - \sum_{e\in E(\gf):\, e^\alpha=x}\frac{(u-1)(e\overline
  e)^\ell}{1 - (u-1)^2(e\overline e)^\ell}\right)^{-1}\in\Bbbk[[u]].$$
(A combinatorial interpretation of these $K_x$ will be given in
Section~\ref{sec:proof}.)
Let $\ell$ be a complete $\Bbbk$-labelling of $\gf$, and define a new
labelling $\ell':E(\gf)\to\Bbbk[[u]]$ by
$$e^{\ell'} = \frac1{1-(u-1)^2(e\overline e)^\ell}e^\ell K_{e^\omega}.$$

Then our main result is the following:
\begin{thm}\label{thm:main}
  With the definitions of $\ell'$ and $K_x$ given above, $\ell'$ is a
  complete labelling and we have in $\Bbbk[[u]]$ the equality
  \begin{equation}\label{eq:a}
    \mF(\ell) = K_\birth\cdot\mG(\ell').
  \end{equation}
\end{thm}

We now explicit the definitions and main result for the length
labelling on a locally finite graph.
\begin{defn}[Path Series]
  The integer-valued series $$G(t)=\sum_{\pi\in[\birth,\death]}t^{|\pi|}\in\N[[t]]$$
  is called the \emdef{path series} of $(\gf,\birth,\death)$.
  The series
  $$F(u,t)=\sum_{\pi\in[\birth,\death]}u^{\bcount(\pi)}t^{|\pi|}\in\N[u][[t]]\subset\N[[u,t]]$$
  is called the \emdef{enriched path series} of $(\gf,\birth,\death)$.
  Its specialization $F(0,t)$ is called the \emdef{proper path series}
  of $(\gf,\birth,\death)$.

  In case $\birth=\death$, we will call $G$ the \emdef{circuit series} of
  $(\gf,\birth)$ and $F$ the \emdef{enriched circuit series} of
  $(\gf,\birth)$. The circuit series is often called the \emph{Green
    function} of the graph $\gf$.
\end{defn}
Note that $F(u,t)$ lies in $\N[u][[t]]$ because the number of bumps on
a path is smaller than its length, so all monomials in the sum have a
$u$-degree smaller than their $t$-degree; hence for any fixed
$t$-degree there are only finitely many monomials with same
$t$-degree, because $\gf$ is locally finite.

Expressed in terms of length labellings, our main theorem then gives
the following result:
\begin{cor}\label{cor:main}
  Suppose $\gf$ is a $d$-regular graph. Then we have
  \begin{equation}\label{eq:b}
    \frac{F(1-u,t)}{1-u^2t^2} =
    \frac{G\left(\frac{t}{1+u(d-u)t^2}\right)}{1+u(d-u)t^2}.
  \end{equation}
\end{cor}

\begin{proof}
  Because $\gf$ is regular, the $K_x$ defined above do not depend on
  $x$ and are all equal to
  $$K = \left(1 - d\frac{(u-1)t^2}{1-(u-1)^2t^2}\right)^{-1} = \frac{1-(1-u)^2t^2}{1+(d-1+u)(1-u)t^2};$$
  thus Theorem~\ref{thm:main} reads
  $$\mF(e\mapsto t) = K\cdot\mG\left(e\mapsto\frac{tK}{1-(1-u)^2t^2}=:\circledast\right).$$
  Now writing $\mF(e\mapsto t)=F(u,t)$ and
  $\mG(e\mapsto\circledast)=G(\circledast)$ completes the proof.
\end{proof}
The special case $u=1$ of this formula appears as an exercise
in~\cite[page~72]{godsil:ac}.

The meaning of the corollary is that, for regular graphs, the richer
two-variable generating series $F(u,t)$ can be recovered from the
simpler $G(t)$. Conversely, $G$ can be recovered from $F$, for
instance because $G(t) = F(1,t)$. Remember it is valid to substitute
$1$ for $u$ in $F$, because for any fixed $t$-degree only finitely
many monomials with that $t$-degree occur in $F$. In fact, much more
is true, as we have the equality
$$G(z) = \frac{1+\frac{\left(1-\sqrt{1-4(1-u)(d-1+u)z^2}\right)^2}{4(1-u)(d-1+u)z^2}}{1-\frac{\left(1-\sqrt{1-4(1-u)(d-1+u)z^2}\right)^2}{4(d-1+u)^2z^2}}F\left(u,\frac{1-\sqrt{1-4(1-u)(d-1+u)z^2}}{2(1-u)(d-1+u)z}\right),$$
or after simplification
\begin{multline}\label{eq:c}
  G(z) = \frac{2(d-1+u)}{d-2+u+(d+u)\sqrt{1-4(1-u)(d-1+u)z^2}}\\
  \cdot F\left(u,\frac{1-\sqrt{1-4(1-u)(d-1+u)z^2}}{2(1-u)(d-1+u)z}\right)
\end{multline}
where both sides are to be understood as power series in $\N[[u,z]]$
that actually reside in $\N[[z]]$. Then for any value (say, in $\C$)
of $u$ we obtain an expression of $G$ in terms of $F(u,-)$. Of
particular interest is the case $u=0$, where (\ref{eq:c}) specializes
to
\begin{equation}\label{eq:c0}
  G(z) = \frac{2(d-1)}{d-2+d\sqrt{1-4(d-1)z^2}}F\left(0,\frac{1-\sqrt{1-4(d-1)z^2}}{2(d-1)z}\right).
\end{equation}
This equation appears in a slightly different form in~\cite{grigorchuk:phd}.

Similarly, we have in $\N[[t,z,z^{-1}]]$ the equality
\begin{equation}\label{eq:d}
  G(z) = \frac{2}{2-d^2tz+dz\sqrt{d^2t^2+4-4t/z}}F\left(1-\frac{dt-\sqrt{d^2t^2+4-4t/z}}{2t},t\right).
\end{equation}
Beware though that (\ref{eq:d}) holds for formal variables $z$ and
$t$; if we were to substitute a real number for $t$, then the
resulting series $G(z)$ would converge absolutely for
$\frac{t}{1+(d-1)t^2}\le z\le t<\rho$, where $\rho$ is the radius of
convergence of $F(1,-)=G$, and in particular not in a neighbourhood of
$0$.

The equalities~(\ref{eq:c}) and~(\ref{eq:d}) are easily derived
from~(\ref{eq:b}) by letting $z=\frac{t}{1+u(d-u)t^2}$ and solving for
$t$ and $u$.

\begin{cor}\label{cor:ratalg}
  In the setting described above:
  \begin{enumerate}
  \item If $\gf$ is finite, then both $F$ and $G$ are rational
    series.
  \item If $G$ is rational, then $F$ is rational too.
  \item $F$ is algebraic if and only if $G$ is algebraic.
  \end{enumerate}
\end{cor}
The proofs are immediate and follow from the explicit form
of~(\ref{eq:b}). The converse of the first statement of the preceding
corollary will be proved in Section~\ref{subs:convc}. The last
statement appears in~\cite{grigorchuk:phd}
and~\cite{grigorchuk-h:groups}.

In the following section we draw some applications to other fields:
group theory, language theory. We give applications of
Theorem~\ref{thm:main} and Corollary~\ref{cor:main} to some examples
of graphs in Section~\ref{sec:examples}, and a derivation of a
``cogrowth formula'' (as that of Subsection~\ref{subs:gptheory}) for a
non-free presentation in Section~\ref{sec:nonfree}.

We give two proofs of the main result in Sections~\ref{sec:newpf}
and~\ref{sec:proof}. The first one, shorter, uses a reduction to
finite graphs and their adjacency matrices. The second one is
combinatorial and uses the inclusion-exclusion principle. 
Using the first proof, we obtain in Section~\ref{sec:matrices} an
extension of a result by Yasutaka Ihara.

Finally in Section~\ref{sec:freeamalg} we show how to compute the
circuit series of a free product of graphs (an analogue of the free
products of groups, \emph{via} their Cayley graph), and in
Section~\ref{sec:direct} do the same for direct products of graphs.

\section{Applications to Other Fields}\label{sec:applic}
The original motivation for Formula~\ref{eq:b} was its implication of
a well-known result in the theory of random walks on discrete groups.

\subsection{Applications to Random Walks on Groups}\label{subs:gptheory}
In this section we will show how $G$ is related to random walks and
$F$ to cogrowth. This will give a generalization of the main
formula~(\ref{eq:grigformula}) to homogeneous spaces $\Pi/\Xi$, where
$\Xi$ does not have to be normal and $\Pi$ is a free product of
infinite-cyclic and order-two groups. For a survey on the topic of
random walks see~\cite{mohar-w:spectra,woess:cogrowth}.

Throughout this subsection we will have $F(t) = F(0,t)$.  We recall
the notion of growth of groups:
\begin{defn}
  Let $\Gamma$ be a group generated by a finite symmetric set $S$.
  For a $\gamma\in\Gamma$ define its \emdef{length}
  $$|\gamma| = \min\{n\in\N:\,\gamma\in S^n\}.$$
  The \emdef{growth series} of $(\Gamma,S)$ is the formal power series
  $$f_{(\Gamma,S)}(t) = \sum_{\gamma\in\Gamma}t^{|\gamma|}\in \N[[t]].$$
  Expanding $f_{(\Gamma,S)}(t)=\sum f_nt^n$, the \emph{growth} of
  $(\Gamma,S)$ is
  $$\alpha(\Gamma,S)=\limsup_{n\to\infty}\sqrt[n]{f_n}$$
  (this supremum-limit is actually a limit and is smaller than $|S|-1$).

  Let $R$ be a subset of $\Gamma$. The \emph{growth series} of $R$
  relative to $(\Gamma,S)$ is the formal power series
  $$f^R_{(\Gamma,S)}(t) = \sum_{\gamma\in R}t^{|\gamma|}\in \N[[t]].$$
  Expanding $f^R_{(\Gamma,S)}(t)=\sum f_nt^n$, define the \emph{growth} of
  $R$ relative to $(\Gamma,S)$ as
  $$\alpha(R;\,\Gamma,S)=\limsup_{n\to\infty}\sqrt[n]{f_n}.$$
  
  If $X$ is a transitive right $\Gamma$-set, the \emdef{simple random
    walk} on $(X,S)$ is the random walk of a point on $X$, having
  probability $1/|S|$ of moving from its current position $x$ to a
  neighbour $x\cdot s$, for all $s\in S$. Fix a point $\birth\in X$,
  and let $p_n$ be the probability that a walk starting at $\birth$
  finish at $\birth$ after $n$ moves (it does not depend on the choice
  of $\birth$). We define the \emdef{spectral radius} of the random
  walk as
  $$\nu(X,S) = \limsup_{n\to\infty}\sqrt[n]{p_n}.$$

  A group $\Pi$ is \emdef{quasi-free} if it is free product of cyclic
  groups of order $2$ and~$\infty$. Equivalently, there exists a
  finite set $S$ and an involution $\overline\cdot:S\to S$ such that,
  as a monoid,
  $$\Pi = \langle S|\,s\overline s=1\;\forall s\in S\rangle.$$
  $\Pi$ is then said to be \emph{quasi-free on $S$}. All quasi-free
  groups on $S$ have the same Cayley graph, which is a regular tree of
  degree $|S|$.

  Every group $\Gamma$ generated by a symmetric set $S$ is a quotient
  of a quasi-free group in the following way: let $\overline\cdot$ be
  an involution on $S$ such that for all $s\in S$ we have the equality
  $\overline s=s^{-1}$ in $\Gamma$. Then $\Gamma$ is a quotient of the
  quasi-free group $\langle S|\,s\overline s=1\;\forall s\in S\rangle$.

  The \emdef{cogrowth series} (respectively \emdef{cogrowth}) of
  $(\Gamma,S)$ is defined as the growth series (respectively growth)
  of $\ker(\pi:\Pi\to\Gamma)$ relative to $(\Pi,S)$, where $\Pi$ is
  a quasi-free group on $S$.
\end{defn}

Associated with a group $\Pi$ generated by a set $S$ and a subgroup
$\Xi$ of $\Pi$, there is a $|S|$-regular graph $\gf$ on which
$\Pi$ acts, called the \emdef{Schreier graph} of $(\Pi,S)$ relative to
$\Xi$. It is given by $\gf=(E,V)$, with
$$V = \Xi\backslash\Pi$$
and
$$E = V\times S,\quad (v,s)^\alpha=v,\quad (v,s)^\omega=vs, \quad
\overline{(v,s)}=(vs,s^{-1});$$
i.e. two cosets $A,B$ are joined by at least one edge if and only if
$AS\supset B$.  (This is the Cayley graph of $(\Pi,S)$ if $\Xi=1$.)
There is a circuit in $\gf$ at every vertex $\Xi v\in\Xi\backslash\Pi$
such that with $s\in v^{-1}\Xi v$ for some $s\in S$; and there is a
multiple edge from $\Xi v$ to $\Xi w$ in $\gf$ if there are $s,t\in
v^{-1}\Xi w$ with $s\neq t\in S$.

\begin{cor}\label{cor:homo}
  Let $\Pi$ be a quasi-free group, presented as a monoid as
  $$\Pi = \langle S|\,s\overline s=1\;\forall s\in S\rangle.$$
  Let $\Xi<\Pi$ be a subgroup of $\Pi$.  Let $\nu=\nu(\Xi\backslash\Pi,S)$
  denote the spectral radius of the simple random walk on
  $\Xi\backslash\Pi$ generated by $S$; and $\alpha=\alpha(\Xi;\,\Pi,S)$
  denote the relative growth of $\Xi$ in $\Pi$. Then we have
  \begin{equation}\label{eq:gchomo}
    \nu = \begin{cases}
      \frac{\sqrt{|S|-1}}{|S|}\left(\frac\alpha{\sqrt{|S|-1}}+\frac{\sqrt{|S|-1}}\alpha\right) & \text{if }\alpha>\sqrt{|S|-1},\\
      \frac{2\sqrt{|S|-1}}{|S|} & \text{if }\alpha\le\sqrt{|S|-1}.
    \end{cases}
  \end{equation}
\end{cor}

\begin{proof}
  Let $\gf$ be the Schreier graph of $(\Pi,S)$ relative to $\Xi$ defined
  above.  Fix the endpoints $\birth=\death=\Xi$, the coset of $1$, and
  give $\gf$ the length labelling. Let $G$ and $F$ be the circuit and
  proper circuit series of $\gf$. In this setting, expressing
  $F(t)=\sum f_n t^n$ and $G(t)=\sum g_n t^n$, we see that $|S|\nu$ is
  the growth rate $\limsup\sqrt[n]{g_n}$ of circuits in $\gf$, and
  $\alpha$ the growth rate $\limsup\sqrt[n]{f_n}$ of proper circuits
  in $\gf$. As both $F$ and $G$ are power series with non-negative
  coefficients, $1/(|S|\nu)$ is the radius of convergence of $G$ and
  $1/\alpha$ the radius of convergence of $F$. Let $d=|S|$ and
  consider the function
  $$(t)\phi = \frac{t}{1+(d-1)t^2}.$$
  This function is strictly increasing for $0\le t<1/\sqrt{d-1}$,
  has a maximum at $t=1/\sqrt{d-1}$ with $(t)\phi=1/(2\sqrt{d-1})$,
  and is strictly decreasing for $t>1/\sqrt{d-1}$.
  
  First we suppose that $\alpha\ge\sqrt{d-1}$, so $\phi$ is monotonously
  increasing on $[0,1/\alpha]$. We set $u=1$ in~(\ref{eq:b})
  and note that, for $t<1$, it says that $F$ has a singularity at $t$
  if and only if $G$ has a singularity at $(t)\phi$. Now as
  $1/\alpha$ is the singularity of $F$ closest to $0$, we conclude by
  monotonicity of $\phi$ that the singularity of $G$ closest to $0$ is
  at $(1/\alpha)\phi$; thus
  $$\frac1{d\nu}=\frac{1/\alpha}{1+(d-1)/\alpha^2}=(1/\alpha)\phi.$$
  
  Suppose now that $\alpha<\sqrt{d-1}$. If $d\nu<2\sqrt{d-1}$,
  the right-hand side of~(\ref{eq:b}) would be bounded for all
  $t\in\R$ while the left-hand side diverges at $t=1$. If
  $d\nu>2\sqrt{d-1}$, there would be a $t\in[0,1/\sqrt{d-1}[$ with
  $(t)\phi=d\nu$; and $F$ would have a singularity at
  $t<1/\alpha$. The only case left is $d\nu=2\sqrt{d-1}$.
\end{proof}
\begin{figure}
\begin{center}
\setlength\unitlength{1pt}
\begin{picture}(300,300)
\put(30,295){\makebox(0,0)[lb]{\smash{\normalsize $\nu$}}}
\put(275,50){\makebox(0,0)[lb]{\smash{\normalsize $\alpha$}}}
\put(-30,-10){\epsfig{file=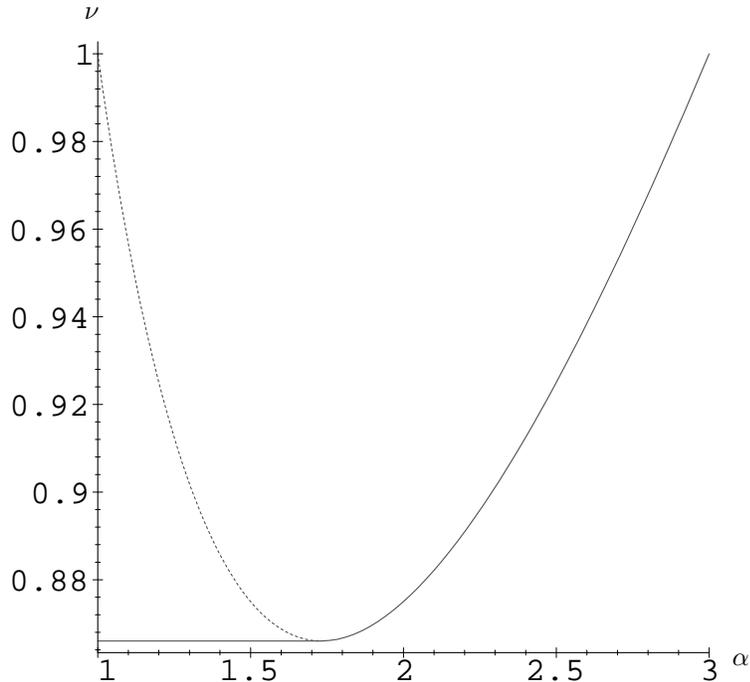}}
\end{picture}
\end{center}
\caption{The function $\alpha\mapsto\nu$ relating cogrowth and spectral
  radius (for $d=4$).}
\end{figure}

\begin{cor}[Grigorchuk~\cite{grigorchuk:rw}]\label{cor:cogrowth}
  Let $\Gamma$ be a group generated by a symmetric finite set $S$, let
  $\nu$ denote the spectral radius of the simple random walk on
  $\Gamma$, and let $\alpha$ denote the cogrowth of $(\Gamma,S)$. Then
  \begin{equation}\label{eq:gc}
    \nu = \begin{cases}
      \frac{\sqrt{|S|-1}}{|S|}\left(\frac\alpha{\sqrt{|S|-1}}+\frac{\sqrt{|S|-1}}\alpha\right) & \text{if }\alpha>\sqrt{|S|-1},\\
      \frac{2\sqrt{|S|-1}}{|S|} & \text{else}.
    \end{cases}
  \end{equation}
\end{cor}
A variety of proofs exist for this result: the original~\cite{grigorchuk:rw}
by Grigorchuk, one by Cohen~\cite{cohen:cogrowth}, an extension by
Northshield to regular graphs~\cite{northshield:cogrowth}, a short proof by
Szwarc~\cite{szwarc:cogrowth} using operator theory, one by
Woess~\cite{woess:cogrowth}, etc.

\begin{proof}
  Present $\Gamma$ as $\Pi/\Xi$, with $\Pi$ a quasi-free group and $\Xi$
  the normal subgroup of $\Pi$ generated by the relators in $\Gamma$,
  and apply Corollary~\ref{cor:homo}.
\end{proof}

We note in passing that if $\alpha<\sqrt{|S|-1}$, then necessarily
$\alpha=0$. Equivalently, we will show that if $\alpha<\sqrt{|S|-1}$,
then $\Xi=1$, so the Cayley graph $\gf$ is a tree. Indeed, suppose
$\gf$ is not a tree, so it contains a circuit $\lambda$ at $\birth$.
As $\gf$ is transitive, there is a translate of $\lambda$ at every
vertex, which we will still write $\lambda$. There are at least
$|S|(|S|-1)^{t-2}(|S|-2)$ paths $p$ of length $t$ in $\gf$ starting at
$\birth$ such that the circuit $p\lambda\overline p$ is proper; thus
$$\alpha \ge \limsup_{t\to\infty}\sqrt[2t+|\lambda|]{|S|(|S|-1)^{t-2}(|S|-2)} = \sqrt{|S|-1}.$$
In fact it is known that $\alpha>\sqrt{|S|-1}$; see~\cite{paschke:norm}.

\subsection{The Series $F$ and $G$ on their Circle of
  Convergence}\label{subs:convc} In this subsection we study the
singularities the series $F$ and $G$ may have on their circle of
convergence. The smallest positive real singularity has a special
importance:
\begin{defn}
  For a series $f(t)$ with positive coefficients, let $\rho(f)$ denote
  its radius of convergence. Then $f$ is \emdef{$\rho(f)$-recurrent}
  if
  $$\lim_{t\to\rho(f)}f(t) = \infty.$$
  Otherwise, it is \emdef{$\rho(f)$-transient}.
\end{defn}
As typical examples, $1/(\rho-t)$ is $\rho$-recurrent, as are all
rational series; $1/\sqrt{\rho-t}$ is not $\rho$-transient.

To study the singularities of $F$ or $G$, we may suppose that
$\birth=\death$; indeed in was shown in~\cite{kesten:rwalks}
and~\cite[Lemma~1]{woess:cogrowthsrw} that the singularities of $F$
and $G$ do not depend on the choice of $\birth$ and $\death$. We make 
that assumption for the remainder of the subsection. We will also
suppose throughout that $\gf$ is $d$-regular, that the radius of
convergence of $F$ is $1/\alpha$ and the radius of convergence of $G$
is $1/(d\nu)=1/\beta$.

\begin{defn}
  Let $\gf$ be a connected graph. A \emdef{proper cycle} in $\gf$ is a
  proper circuit $(\pi_1,\dots,\pi_n)$ such that
  $\overline{\pi_1}\neq\pi_n$. The \emdef{proper period} $p$ and
  \emdef{strong proper period} $p_s$ are defined as follows:
  \begin{gather*}
    p = \gcd\{n|\,\text{there exists a proper cycle $\pi$ in
      $\gf$ with }|\pi|=n\},\\
    p_s=\gcd\{n|\,\forall x\in V(\gf)\text{ there exists a
      proper cycle $\pi$ in $\ball(x,n)$ with }|\pi|=n\},
  \end{gather*}
  where by convention the $\gcd$ of the empty set is $\infty$.  The
  graph $\gf$ is \emdef{strongly properly periodic} if $p=p_s$.
  
  The \emdef{period} $q$ and \emdef{strong period} $q_s$ of $\gf$ are
  defined analogously with `proper cycle' replaced by `circuit'. $\gf$ is
  \emdef{strongly periodic} if $q=q_s$.
\end{defn}
  
\begin{thm}[Cartwright~\cite{cartwright:sing}]\label{thm:cart}
  Let $\gf$ have proper period $p$ and strong proper period
  $p_s$. Then the singularities of $F$ on its circle of convergence
  are among the
  $$\frac{e^{2i\pi k/p_s}}{\alpha},\quad k=1,\dots,p_s.$$
  If moreover $\gf$ is strongly properly periodic, the singularities
  of $F$ on its circle of convergence are precisely these numbers.
  
  Let $\gf$ have period $q$ and strong period $q_s$. Then the
  singularities of $G$ on its circle of convergence are among the
  $$\frac{e^{2i\pi k/q_s}}{\beta},\quad k=1,\dots,q_s.$$
  If moreover $\gf$ is strongly periodic, the singularities of $G$ on
  its circle of convergence are precisely these numbers.
\end{thm}

If $\gf$ is connected and non-trivial, there is a path of even length
at every vertex (a sequence of bumps, for instance). All graphs are
then either $2$-periodic (if they are bipartite) or $1$-periodic. If
there is a constant $N$ such that for all $x\in V(\gf)$ there is at
$x$ a circuit of odd length bounded by $N$, then $\gf$ is strongly
$1$-periodic; otherwise $\gf$ is strongly $2$-periodic. The
singularities of $G$ on its circle of convergence are then at
$1/\beta$, and also at $-1/\beta$ if $\gf$ is strongly periodic with
period $2$.

If $\gf$ is not strongly periodic, there may be one or two
singularities on $G$'s circle of convergence; consider for instance
the $4$-regular tree, and at a vertex $\birth$ delete
two or three edges replacing them by loops. The resulting graphs
$\gf_2$ and $\gf_3$ are still $4$-regular and their circuit series, as
computed in~(\ref{eq:treeloop}), are respectively
\begin{equation}\label{ex:4reg3}
  G_2(t) = \frac3{2-6t+\sqrt{1-12t^2}}\quad\text{and}\quad G_3(t) = \frac6{5-18t+\sqrt{1-12t^2}}.
\end{equation}
$G_2$ has singularities at $\pm1/\sqrt{12}$ on its circle of
convergence, while $G_3$ has only $2/7$ as singularity on its circle
of convergence.

Following the proof of Corollary~\ref{cor:homo} above, we see
that if $\beta<d$ the singularities of $F$ on its circle of convergence
are in bijection with those of $G$, so are at $1/\alpha$ and possibly
$-1/\alpha$, if $\gf$ is strongly two-periodic. If $\beta=d$, though,
$\gf$ can have any strong proper period; consider for example the
cycles on length $k$ studied in Section~\ref{subs:cycles}: they are
strongly properly $k$-periodic.

This simple result shows how $\gf$ can be approximated by finite
subgraphs.
\begin{lem}\label{lem:limit}
  Let $\gf$ be a graph and $x$, $y$ two vertices in $\gf$.  Let
  $\mG_{x,y}$ and $\mF_{x,y}$ be the path series and enriched path
  series respectively from $x$ to $y$ in $\gf$, and let $\mG_{x,y}^n$
  and~$\mF_{x,y}^n$ be the path series and enriched path series
  respectively from $x$ to $y$ in the ball $\ball(x,n)$ (they are $0$
  if $y\not\in\ball(x,n)$). Then
  $$\lim_{n\to\infty}\mG_{x,y}^n = \mG_{x,y},\qquad\lim_{n\to\infty}\mF_{x,y}^n = \mF_{x,y}.$$
\end{lem}
\begin{proof}
  Recall that $\lim\mG_{x,y}^n=\mG_{x,y}$ means that both terms are
  sums of paths, say $A_n$ and $A$, such that the minimal length of
  paths in the symmetric difference $A_n\triangle A$ tends to
  infinity. Now the difference between $\mG_{x,y}^n$ and $\mG_{x,y}$
  consists only of paths in $\gf$ that exit $\ball(x,n)$, and thus
  have length at least $2n-\dist(x,y)\to\infty$. The same argument
  holds for $\mF$.
\end{proof}

\begin{defn}
  The graph $\gf$ is \emdef{quasi-transitive} if $\Aut(\gf)$ acts with
  finitely many orbits.
\end{defn}

\begin{lem}\label{lem:gfcard} Let $\gf$ be a regular quasi-transitive
  connected graph with distinguished vertex $*$, and let $f_n$ and
  $g_n$ denote respectively the number of proper circuits and circuits
  at $*$ of length $n$. Then 
  $$\limsup_{n\to\infty}g_n/\beta^n = \limsup_{n\to\infty}f_n/\alpha^n = \begin{cases}1/|\gf|&\text{if $\gf$ is finite and has odd circuits;}\\
      2/|\gf|&\text{if $\gf$ is finite and has only even circuits;}\\
      0&\text{if $\gf$ is infinite.}
    \end{cases}$$
\end{lem}
\begin{proof}
  If $\gf$ is finite, then $\beta=d$, the degree of $\gf$; after a large
  even number of steps, a random walk starting at $\birth$ will be
  uniformly distributed over $\gf$ (or over the vertices at even
  distance of $\birth$, in case all circuits have even length). A long
  enough walk then has probability $1/|\gf|$ (or $2/|\gf|$ if all
  circuits have even length) of being a circuit.
  
  If $\gf$ is infinite, we consider two cases. If $G(1/\beta)<\infty$,
  i.e.\ $G$ is $1/\beta$-transient, the general term $g_n/\beta^n$ of
  the series $G(1/\beta)$ tends to $0$. If $G$ is $1/\beta$-recurrent,
  then, as $\gf$ is quasi-transitive, $\beta=d$
  by~\cite[Theorem~7.7]{woess:rw}. We then approximate $\gf$ by the
  sequence of its balls of radius $R$:
  $$\lim_{n\to\infty}\frac{g_n}{\beta^n}=\lim_{R,n\to\infty}\frac{g_{R,n}}{d^n}=\lim_{R\to\infty}\frac{(1\text{ or }2)}{|\ball(\birth,R)|}=0,$$
  where we expand the circuit series of $\ball(\birth,R)$ as $\sum g_{R,n}t^n$.

  The same proof holds for the $f_n$. Its particular case where $\gf$
  is a Cayley graph appears in~\cite{woess:cogrowthsrw}.
\end{proof}
Note that if $\gf$ is not quasi-transitive, a somewhat weaker result
holds~\cite[\S 7.1]{kitchens:symbdyn}: if $\gf$ is transient or
null-recurrent then the common limsup is $0$. If $\gf$ is
positive-recurrent then the limsups are normalized coefficients of
$\gf$'s Perron-Frobenius eigenvector.  Lemma~\ref{lem:gfcard} is not
true for arbitrary $d$-regular graphs: consider for instance the graph
$\gf_3$ described above. Its circuit series $G_3$, given
in~(\ref{ex:4reg3}), has radius of convergence $1/\beta=2/7$, and one
easily checks that all its coefficients $g_n$ satisfy
$g_n/\beta^n\ge1/2$.

We obtain the following characterization of rational series:
\begin{thm}
  For regular quasi-transitive connected graphs, the following are
  equivalent:
  \begin{enumerate}
  \item $\gf$ is finite;\label{thm:rat:a}
  \item $G(t)$ is a rational function of $t$;\label{thm:rat:b}
  \item $F(t)$ is a rational function of $t$, and $\gf$ is not an
    infinite tree.\label{thm:rat:c}
  \end{enumerate}
\end{thm}
\begin{proof}
  By Corollary~\ref{cor:ratalg}, statement \ref{thm:rat:a} implies the
  other two. By Corollary~\ref{cor:main}, and a computation on trees
  done in section~\ref{subs:trees} to deal with the case $F(t)=1$,
  statement \ref{thm:rat:b} implies~\ref{thm:rat:c}.  It remains to
  show that statement~\ref{thm:rat:c} implies~\ref{thm:rat:a}.
  
  Assume that $F(t)=\sum f_nt^n$ is rational, not equal to $1$. As the
  $f_n$ are positive, $F$ has a pole, of multiplicity $m$, at
  $1/\alpha$. There is then a constant $a>0$ such that
  $f_n>a\binom{n}{m-1}\alpha^n$ for infinitely many values of
  $n$~\cite[page~341]{graham-k-p:concmath}. It follows by
  Lemma~\ref{lem:gfcard} that $m=1$ and
  the graph $\gf$ is finite, of cardinality at most $1/a$.
\end{proof}
It is not known whether the same holds for regular, or even arbitrary
connected graphs. Certainly an altogether different proof would be
needed.

\subsection{Application to Languages}\label{sec:lang}
Let $S$ be a finite set of cardinality $d$ and let $\overline\cdot$ be
an involution on $S$. A \emdef{word} is an element $w$ of the free
monoid $S^*$. A \emdef{language} is a set $L$ of words. The language
$L$ is called \emdef{saturated} if for any $u, v\in S^*$ and $s\in S$
we have
$$uv\in L \Longleftrightarrow us\overline sv\in L;$$
that is to say, $L$ is stable under insertion and deletion of subwords
of the form $s\overline s$. The language $L$ is called \emdef{desiccated} if
no word in $L$ contains a subword of the form $s\overline s$. Given a
language $L$ we may naturally construct its \emph{saturation} $\langle
L\rangle$, the smallest saturated language containing $L$, and its
\emph{desiccation} $\widehat L$, the largest desiccated language
contained in $L$.

Let $\Sigma$ be the monoid defined by generators $S$ and relations
$s\overline s=1$ for all $s\in S$:
\begin{equation}\label{eq:defSigma}
  \Sigma = \langle S |\, s\overline s=1\;\forall s\in S\rangle.
\end{equation}
This is a free product of free groups and order-two groups; if
$\overline\cdot$ is fixed-point-free, $\Sigma$ is a free group.
Write $\phi$ for the canonical projection from $S^*$ to $\Sigma$.
Let $\Bbbk=\Z[\Sigma]$ be its monoid ring. Then given a language
$L\subset S^*$ we may define its \emph{growth series}
$\Theta(L)$ as
$$\Theta(L) = \sum_{w\in L} w^\phi t^{|w|}\in \Bbbk[[t]].$$

This notion of growth series with coefficients was introduced by
Fabrice Liardet in his doctoral thesis~\cite{liardet:thesis}, where he
studied \emph{complete growth functions} of groups.

\begin{thm}
  For any language $L$ there holds
  \begin{equation}\label{eq:lang}
    \frac{\Theta(\widehat L)(t)}{1-t^2} =
    \frac{\Theta(\langle L\rangle)\left(\frac{t}{1+(d-1)t^2}\right)}{1+(d-1)t^2},
  \end{equation}
  where $d = |S|$.
\end{thm}

\begin{proof}
  For any language there exists a unique minimal (possibly infinite)
  automaton recognising it~(\cite[\S III.5]{eilenberg:a} is a good
  reference). Let $\gf$ be the minimal automaton recognising $\langle
  L\rangle$. Recall that this is a graph with an initial vertex
  $\birth$, a set of terminal vertices $T$ and a labelling
  $\ell':E(\gf)\to S$ of the graph's edges such that the number of
  paths labelled $w$, starting at $\birth$ and ending at a $\tau\in T$
  is $1$ if $w\in L$ and $0$ otherwise.  Extend the labelling $\ell'$
  to a labelling $\ell:E(\gf)\to \Bbbk[[t]]$ by
  $$e^\ell = t\cdot(e^{\ell'})^\phi.$$

  Because $\langle L\rangle$ is saturated, and $\gf$ is minimal,
  $(\overline e)^\ell=\overline{e^\ell}$; then $\widehat L$ is
  the set of labels on proper paths from $\birth$ to some $\tau\in T$.
  Choosing in turn all $\tau\in T$ as $\death$, we obtain growth series
  $F_\tau, G_\tau$ counting the formal sum of paths and proper paths from
  $\birth$ to $\tau$. It then suffices to write
  $$\frac{\Theta(\widehat L)(t)}{1-t^2} = \frac{\sum_{\tau\in T}F_\tau(t)}{1-t^2}
  = \frac{\sum_{\tau\in T}G_\tau\left(\frac{t}{1+(d-1)t^2}\right)}{1+(d-1)t^2} = \frac{\Theta(\langle L\rangle)\left(\frac{t}{1+(d-1)t^2}\right)}{1+(d-1)t^2}.$$
\end{proof}

The following result is well-known:
\begin{thm}[M\"uller\&Schupp~\cite{muller-s:context-free,muller-s:context-free-2}]
  Let $\Gamma$ be a finitely generated group, presented as a quotient
  $\Sigma/\Xi$ with $\Sigma$ as in~(\ref{eq:defSigma}). Then
  $\Theta(\Xi)$ is an algebraic series (i.e.\ satisfies a polynomial
  equation) if and only if $\Sigma/\Xi$ is virtually free (i.e.\ has a
  normal subgroup of finite index that is free).
\end{thm}
It is not known whether there exists a non-virtually-free graph whose
circuit series (as defined in Corollary~\ref{cor:main}) is algebraic.

\section{First Proof of Theorem~\ref{thm:main}}\label{sec:newpf}
We will now prove Theorem~\ref{thm:main} using linear algebra.  We
first assume the graph has a finite number of vertices; for the
computations refer to $\Bbbk$-matrices and $\Bbbk[[u]]$-matrices
indexed by the graph's vertices. This proof is hinted at in Godsil's
book as an exercise~\cite[page~72]{godsil:ac}; it was suggested to the
author by Gilles Robert.

For all pairs of vertices $x,y\in V(\gf)$ let
$$\mG_{x,y}(\ell)=\sum_{\pi\in[x,y]}\pi^\ell,\qquad\mF_{x,y}(\ell)=\sum_{\pi\in[x,y]}u^{\bcount(\pi)}\pi^\ell$$
be the path and enriched path series from $x$ to $y$; for ease of
notation we will leave out the labelling $\ell$ if it is obvious from
context.  Let $\delta_{x,y}$ denote the Kronecker delta, equal to $1$
if $x=y$ and $0$ otherwise. For any $v\in\Bbbk$, let $[v]_x^y$ denote
the $V(\gf)\times V(\gf)$ matrix with zeroes everywhere except at
$(x,y)$, where it has value $v$. Then
$$\mG_{x,y} = \delta_{x,y} + \sum_{e\in E(\gf): e^\alpha = x} e^\ell\mG_{e^\omega,y}$$
so that if
$$A = \sum_{e\in E(\gf)} [e^\ell]_{e^\alpha}^{e^\omega}$$
be the adjacency matrix of $\gf$, with labellings, then we have
$$(\mG_{x,y})_{x,y\in V(\gf)} = \frac1{1 - A},$$
an equation holding between $V(\gf)\times V(\gf)$ matrices over $\Bbbk$.

Similarly, letting $\mF_{x,e,y}$ count the paths from $x$ to $y$ that
start with the edge $e$,
\begin{align*}
\mF_{x,y} &= \delta_{x,y} + \sum_{e\in E(\gf): e^\alpha = x}\mF_{x,e,y},\\
\mF_{x,e,y} &= e^\ell\left(\mF_{e^\omega,y} + (u-1)\mF_{e^\omega,\overline e,y}\right),\\
\mF_{e^\omega,\overline e,y} &= \overline e^\ell\left(\mF_{x,y} + (u-1)\mF_{x,e,y}\right);
\end{align*}
these last two lines solve to
$$\mF_{x,e,y} = \left(1 - (u-1)^2(e\overline e)^\ell\right)^{-1}\left(e^\ell\mF_{e^\omega,y} + (u-1)(e\overline e)^\ell\mF_{x,y}\right),$$
which we insert in the first line to obtain
$$K_x^{-1}\mF_{x,y} = \delta_{x,y} + \sum_{e\in E(\gf): e^\alpha = x}\frac{e^\ell}{1-(u-1)^2(e\overline e)^\ell}K_{e^\omega}\cdot K_{e^\omega}^{-1}\mF_{e^\omega,y}.$$
Thus if we let
\begin{equation}\label{eq:defb}
  e^{\ell'}=\frac{e^\ell}{1-(u-1)^2(e\overline e)^\ell}K_{e^\omega},\qquad
  A' = \sum_{e\in E(\gf)}[e^{\ell'}]_{e^\alpha}^{e^\omega},
\end{equation}
we obtain
\begin{equation}\label{eq:feqb}
  (K_x^{-1}\mF_{x,y})_{x,y\in V(\gf)} = \frac1{1-A'}
\end{equation}
and the proof is finished in the case that $\gf$ is finite, because
the matrix $A'$ is precisely that obtained from $A$ by substituting
$\ell'$ for $\ell$.

If $\gf$ has infinitely many vertices, we approximate it, thanks to
Lemma~\ref{lem:limit}, by finite graphs. Denote by
$\mF_{\birth,\death}^n(\ell)$ and $\mG_{\birth,\death}^n(\ell')$ the
enriched path series and path series respectively in $\ball(\birth,n)$, and
write
$$K_\birth\cdot\mF(\ell) = \lim_{n\to\infty}\mF_{\birth,\death}^n(\ell) =
\lim_{n\to\infty}\mG_{\birth,\death}^n(\ell') = \mG(\ell')$$
to complete the proof.

\section{Graphs and Matrices}\label{sec:matrices}
Graphs can be studied through their \emph{adjacency} and
\emph{incidence} matrices. We give here the relevant definitions and
obtain an extension of a theorem by Hyman Bass~\cite{bass:ihara} on
the Ihara-Selberg zeta function. We will use power series with
coefficients in a matrix ring, and fractional expressions in matrices;
by convention, we understand `$X/Y$' as `$X\cdot Y^{-1}$'.

\begin{defn}
  Let $\gf$ be a finite graph.
  The \emdef{edge-adjacency} and \emdef{inversion} matrices of $\gf$,
  respectively $B$ and $J$, are $E(\gf)\times E(\gf)$ matrices over
  $\Z$ defined by
  $$B_{e,f}=\begin{cases}1&\text{if }e^\omega=f^\alpha\\0&\text{else},\end{cases}\qquad
  J_{e,f}=\begin{cases}1&\text{if }\overline e=f\\0&\text{else}.\end{cases}$$
  The \emdef{vertex-adjacency} and \emdef{degree} matrices of $\gf$,
  respectively $A$ and $D$, are $V(\gf)\times V(\gf)$ matrices over
  $\Z$ defined by
  $$A_{v,w}=|\{e\in E(\gf)|\,e^\alpha=v\text{ and }e^\omega=w\}|,\qquad
  D_{v,w}=\begin{cases}\deg(v)&\text{if }v=w,\\ 0&\text{else}.\end{cases}$$
  
  A \emdef{cycle} is the equivalence class of a circuit under cyclic
  permutation of its edges. A \emdef{proper cycle} is a cycle all of
  whose representatives are proper circuits. A cycle is
  \emdef{primitive} if none of its representatives can be written as
  $\pi^k$ for some $k\ge2$. The \emdef{cyclic bump count}
  $\cbcount(\pi)$ of a circuit $\pi=(\pi_1,\dots,\pi_n)$ is
  $$\cbcount(\pi)=|\{i=1,\dots,n|\,\pi_i=\overline{\pi_{i+1}}\}|,$$
  where the edge $\pi_{n+1}$ is understood to be $\pi_1$.
\end{defn}

The matrices given above are related to paths in $\gf$ as follows:
Consider first the matrix
$$M = \idmatrix - (B-(1-u)J)t.$$
Then the $(e,f)$ coefficient of $M^{-1}$ is precisely
$$\sum_{\pi:\pi_1=e,\pi^\omega=f^\alpha}u^{\bcount(\pi f)}t^{|\pi|}.$$
This is clear because the series expansion of $M^{-1}$ is the sum of
sequences of $(B-J)t$ (contributing edges with no bump) and $Jut$
(contributing edges with bumps), with an extra factor of $u$ in case
the path ends in $\overline f$. Likewise, consider the matrix
$$P = \idmatrix - At + (1-u)(D-(1-u)\idmatrix)t^2,$$
whose $(v,w)$-coefficient counts paths from vertex $v$ to vertex $w$.

We now state and prove an extension of a theorem by Bass~\cite{bass:ihara,foata-z:bass,northshield:ihara}:
\begin{thm}\label{thm:bass}
  Let $\mathcal C$ be a set of representatives of primitive cycles in
  $\gf$, and form the \emph{zeta function} of $\gf$
  $$\zeta(u,t) = \prod_{\gamma\in\mathcal C}\frac1{1-u^{\cbcount(\gamma)}t^{|\gamma|}}.$$
  (The choice of representatives does not change the zeta function.)
  Then $\zeta^{-1}$ is a polynomial in $u$ and $t$ and can be expressed as
  \begin{align}
    \frac1{\zeta(u,t)} &= \det M\label{eq:bass1}\\
    &= (1+(1-u)t)^{n}(1-(1-u)^2t^2)^{m-|V(\gf)|}\det P\label{eq:bass2},
  \end{align}
  where
  $$n=|\{e\in E(\gf)|\,e=\overline e\}|,\qquad2m=|\{e\in E(\gf)|\,e\neq\overline e\}|.$$
\end{thm}

The special case $u=n=0$ of this result was stated and proved in the
given sources. We will prove the general statement, using a result of
Shimson Amitsur:
\begin{thm}[Amitsur~\cite{amitsur:sum,reutenauer-s:sum}]
  Let $X_1,\dots,X_k$ be square matrices of the same dimension over an
  arbitrary ring.  Let $S$ contain one representative up to cyclic
  permutation of words over the alphabet $\{1,\dots,k\}$ that are
  \emph{primitive}, i.e.\ such that none of their cyclic permutations
  are proper powers of a word ($S$ is infinite as soon as $k>1$). For
  $p=i_1\dots i_n\in S$ set $X_p=X_{i_1}\dots X_{i_n}$. Then
  $$\det(\idmatrix-(X_1+\dots+X_k)t) = \prod_{p\in S}\det(\idmatrix-X_pt^{|p|}),$$
  considered as an equality of power series in $t$ over the matrix ring.
\end{thm}

The equality~(\ref{eq:bass1}) then follows; indeed, for all edges
$e\in E(\gf)$ let $X_e$ be the $E(\gf)\times E(\gf)$ matrix whose
$e$-th row is the $e$-th row of $B-(1-u)J$, and whose other rows are
$0$. Then clearly $\idmatrix-\sum_{e\in E(\gf)}X_et = M$ and, for any
sequence of edges $\pi$,
$$\det(\idmatrix-X_\pi t^{|\pi|}) = \begin{cases}1-u^{\cbcount(\pi)}t^{|\pi|}&\text{
    if $\pi$ is a circuit,}\\1&\text{else,}\end{cases}$$
so equality of $\zeta(u,t)$ and $\det M$ follows from Amitsur's
Theorem.

To prove~(\ref{eq:bass2}), we consider block matrices of dimension
$|E(\gf)|+|V(\gf)|$. We write the adjacency operator and its adjoint
$$T_{e,v}=\begin{cases}1 & \text{ if }e^\alpha=v,\\ 0 &\text{ else},\end{cases},
\qquad T_{v,e}^*=\begin{cases}1 & \text{ if }e^\alpha=v,\\ 0 &\text{ else},\end{cases}$$
so that
\[B=JTT^*,\qquad D=T^*T,\qquad A=T^*JT,\]
whence
\begin{align*}
  (1-(1-u)^2t^2)^{|V(\gf)|}\det M &= \begin{vmatrix}M & *\\ 0 &
    1-(1-u)^2t^2\end{vmatrix}\\
  &= \begin{vmatrix}1+(1-u)tJ & -JTt\\ -T^* & \idmatrix\end{vmatrix}\cdot
  \begin{vmatrix}\idmatrix & 0\\ T^* & 1-(1-u)^2t^2\end{vmatrix}\\
  &= \begin{vmatrix}\idmatrix & 0\\ T^*(1-(1-u)tJ) & 1-(1-u)^2t^2\end{vmatrix}\cdot
  \begin{vmatrix}1+(1-u)tJ & -JTt\\ -T^* & \idmatrix\end{vmatrix}\\
  &= \begin{vmatrix}1+(1-u)tJ & *\\ 0 & P\end{vmatrix}\\
  &= (1+(1-u)t)^n(1-(1-u)^2t^2)^m\det P.
\end{align*}

\section{Second Proof of Theorem~\ref{thm:main}}\label{sec:proof}
Let $P=[\birth,\death]$ be the set of paths in $\gf$ from $\birth$ to
$\death$. As the principle of inclusion-exclusion will be
applied~\cite{wilf:gf}, it will be helpful to compute in $\Pi=\Z[[P]]$,
the $\Z$-module of functions from the set of paths to $\Z$. We embed
subsets of $P$ in $\Pi$ by mapping a subset to its characteristic
function:
$$P\supset A\mapsto \chi_A,\quad\text{with }
(\pi)\chi_A=\begin{cases}1&\text{if }\pi\in A,\\ 0 &
  \text{otherwise}.\end{cases}$$
Let $\mathcal B$ be the subset of bounded non-negative elements of
$\Pi$ (i.e. functions $f$ such that there is a constant $N$ with
$0\le(\pi)f<N$ for all paths $\pi$). If $\ell$ is a complete labelling
of $\gf$, there is an induced labelling $\ell_*:\mathcal B\to\Bbbk$ given by
$$(f)\ell_* = \sum_{\pi\in P}(\pi)f\pi^\ell.$$
Note that the sum, although infinite, defines an element of $\Bbbk$
due to the fact $\ell$ is complete.

\begin{defn}[Bump Scheme]
  Let $e\in E(\gf)$ and $v\in V(\gf)$. A \emdef{squiggle along $e$} is
  a sequence $(e,\overline e,\dots,e,\overline e)$. A \emdef{squiggle
    at $v$} is a squiggle along $e$ for some edge $e$ such that
  $e^\alpha=v$.
  
  Let $\pi=(v_0,e_1,\dots,e_n,v_n)$ be a path of length $n$ in
  $\gf$. A \emdef{bump scheme} for $\pi$ is a pair
  $B=\big((\beta_0,\dots,\beta_n),(\gamma_1,\dots,\gamma_n)\big)$, with
  \begin{itemize}
  \item for all $i\in\{0,\dots,n\}$, a finite (possibly empty) sequence
    $\beta_i=(\beta_{i,1},\dots,\beta_{i,t_i})$ of squiggles at $v_i$;
  \item for all $i\in\{1,\dots,n\}$, a squiggle $\gamma_i$ along $e_i$.
  \end{itemize}
  The \emdef{weight} $|B|$ of the bump scheme $B$ is defined as
  $$|B| = \sum_{i=0}^n\sum_{j=1}^{t_i} (|\beta_{i,j}|-1) +
  \sum_{i=1}^n |\gamma_i|.$$
\end{defn}

Given a path $\pi$ and a bump scheme $B=(\beta,\gamma)$ for $\pi$,
we obtain a new path $\pi\vee B\in P$, by setting
$$\pi\vee B = \beta_{0,1}\cdots\beta_{0,t_0}\gamma_1e_1\beta_{1,1}
\cdots\gamma_ne_n\beta_{n,1}\cdots\beta_{n,t_n},$$
where the product denotes concatenation.


We now define a linear map $\phi:\Pi\to\Pi[[u]]$ by setting,
for $f\in\Pi$ and $\pi\in P$,
$$(\pi)\big((f)\phi\big) = \sum_{(\rho,B):\,\rho\vee B=\pi}(u-1)^{|B|}(\rho)f,$$
where the sum ranges over all pairs $(\rho,B)$ where $\rho\in P$ and
$B$ is a bump scheme for $\rho$ such that $\rho\vee B=\pi$. Note that
the sum is finite because the edges of $\rho$ and of $B$ form subsets
of those of $\pi$.

\begin{lem}
  For any path $\pi$ we have
  \begin{equation}\label{eq:lempf2}
    (\pi)\big((\chi_P)\phi\big) = u^{\bcount(\pi)}.
  \end{equation}
\end{lem}
\begin{proof}
  Say $\pi=(\pi_1,\dots,\pi_n)$ has $m\geq0$ bumps, at indices
  $b_1,\dots,b_m$ so that $\pi_{b_i} = \overline{\pi_{b_{i+1}}}$. We
  will show the evaluation at $\pi$ of the left-hand side
  of~(\ref{eq:lempf2}) yields $u^m$.
  
  We claim there is a bijection between the subsets $C$ of
  $\{1,\dots,m\}$ and the pairs $(\rho_C,B_C)$ where $\rho_C$ is a
  path and $B_C$ is a bump scheme for $\rho_C$ with $\pi = \rho_C\vee
  B_C$; and further $|B_C|=|C|$.
  
  First, take a $\rho$ and a $B=(\beta,\gamma)$ such that $\rho\vee
  B=\pi$. The path $\rho\vee B$ is obtained by shuffling together the
  edges of $\rho$ and $B$, and this partitions the edges of $\pi$ in
  two classes, namely \emph{(i)} those coming from $\rho$ and
  \emph{(ii)} those coming from $\beta$ and $\gamma$.  Let
  $C\subset\{1,\dots,m\}$ be the indices of the bumps $b_i$ in $\pi$
  coming from $B$, i.e.\ such that $\pi_{b_i}$ and $\pi_{b_{i+1}}$
  belong to the class \emph{(ii)}. One direction of the bijection is
  then $(\rho,B)\mapsto C$.
  
  Conversely, given a subset $C$ consider the set $D=\{b_i|\,i\in
  C\}$. Split it in maximal-length runs of consecutive integers
  $D=D_1\sqcup\dots\sqcup D_t$. For all runs $D_i$ we do the
  following: to $D_i=\{j,j+1,\dots,j+2k-1\}$ of even cardinality we
  associate a squiggle $\gamma_j$ of length $2k$ along $\pi_j$; to
  $D_i=\{j,j+1,\dots,j+2k-2\}$ of odd cardinality we associate a
  squiggle $\beta_{j,l}$ of length $2k$ at $v_{j-1}$; then we delete
  in $\pi$ the edges $\pi_j,\dots,\pi_{j+2k-1}$. This process
  constructs a bump scheme $B=(\beta,\gamma)$ while pruning edges of
  $\pi$, giving a path $\gamma$ with $\gamma\vee B=\pi$. These two
  constructions are inverses, proving the claimed bijection.

  It now follows that
  $$(\pi)(\chi_P)\phi = \sum_{C\in\{1,\dots,m\}}(u-1)^{|B_C|} =
  \sum_{r=0}^m(u-1)^r\binom mr = u^m.$$
\end{proof}

Let $\ell':E(\gf)\to\Bbbk[[u]]$ be defined by
$$e^{\ell'} = \frac1{1-(e\overline e)^\ell(1-u)^2}e^\ell K_{e^\omega}.$$
We prove Theorem~\ref{thm:main} by noting that
$\mG(\ell)=(\chi_P)\ell_*$, that $\mF(\ell)=(\chi_P\phi)\ell_*$,
and that for any $f\in\Pi$ we have $(f\phi)\ell_* = K_\birth(f)\ell'_*$.
To prove this last equality, take a path $\pi=(\pi_1,\dots,\pi_n)$ on
vertices $v_0,\dots,v_n$. Then
$$(\chi_{\{\pi\}}\phi)\ell_* = \sum_B(u-1)^{|B|}(\pi\vee B)^\ell,$$
where the sum ranges over all bump schemes for $\pi$, and
$$K_\birth\pi^{\ell'}=K_{v_0}\frac1{1-(u-1)^2(\pi_1\overline\pi_1)^\ell}\pi_1^\ell K_{v_1}\cdots\frac1{1-(u-1)^2(\pi_n\overline\pi_n)^\ell}\pi_n^\ell K_{v_n}.$$
It is clear these last two lines are equal; for the power series
expansion of the $K_{v_i}$ correspond to all the possible squiggle
sequences at $v_i$, and the power series expansion of the
$1/(1-(u-1)^2(\pi_i\overline\pi_i)^\ell)$ correspond to all possible
squiggles along $\pi_i$.

\section{Examples}\label{sec:examples}
We give here examples of regular graphs and when possible compute
independently the series $F$ and $G$. In some cases it will be easier
to compute $F$, while in others it will be simpler to compute $G$
first. In all cases, once one of $F$ and $G$ has been computed, the
other one can be obtained from Corollary~\ref{cor:main}.

In all the examples the graphs are vertex transitive, so the choice of
$\birth$ is unimportant. To simplify the computations we choose
$\death=\birth$ and the length labelling.

\subsection{Complete Graphs}\label{subs:complete}
Let $\gf=K_v$, the complete graph on $v\ge3$ vertices. Its degree is $d=v-1$.
To compute $F$ and $G$, choose three distinct vertices $\birth,\$,\#$
(the choice is unimportant as $K_v$ is three-transitive). Define growth series
\begin{description}
\item[$F(u,t)$] the growth series of circuits based at $\birth$;
\item[$F'(u,t)$] the growth series of paths $\pi$ from $\$$ to $\birth$ with $\pi_1^\omega=\#$;
\item[$F''(u,t)$] the growth series of paths $\pi$ from $\$$ to $\birth$ with $\pi_1^\omega=\birth$.
\end{description}

Then
\begin{align*}
F &= 1 + (v-1)t \left[(v-2)F' + uF''\right],\\
F' &= t\left[F'' + (v-3+u)F'\right],\\
F'' &= t\left[1 + (F-1)\frac{v-2+u}{v-1}\right];
\end{align*}
Indeed the first line states that a circuit at $\birth$ is either the trivial
circuit at $\birth$, or a choice of one of $v-1$ edges to another point
(call it $\$$), followed by a path from $\$$ to $\birth$; this path can
first go to any vertex of the $v-2$ vertices (say ($\#$) different
from $\birth$ and $\$$, and thus contribute $F'$, or go back to $\birth$ and
contribute $F''$ and a bump.

The second equation says that a path from $\$$ to $\birth$ starting by
going to $\#$ can either continue to $\birth$, contributing $F'$, go to any
of the $v-3$ other vertices contributing $F'$, or come back to
$\$$, contributing $F'$ and a bump.

The third line says that a path from $\$$ to $\birth$ starting by going to
$\birth$ continues as a circuit at $\birth$; but if the circuit is non-trivial, then
one out of $v-1$ times a bump will be contributed.

Solving the system, we obtain
$$F(u,t) = \frac{1+(1-u)t}{1-(v-2+u)t}\cdot\frac{1-(v-2)t+(1-u)(v-2+u)t^2}{1+t+(1-u)(v-2+u)t^2}.$$

We then compute
\begin{gather*}
G(t) = F(1,t) = \frac{1-(v-2)t}{(1+t)(1-(v-1)t)},\\
F(0,t) = \frac{(1+t)(1-(v-2)t+(v-2)t^2)}{(1-(v-2)t)(1+t+(v-2)t^2)}.
\end{gather*}

\subsection{Cycles}\label{subs:cycles}
Let $\gf=C_k$, the cycle on $k$ vertices. Here, as there are $2$ proper
circuits of length $n$ for all $n$ multiples of $k$, we have
$$F(0,t) = \frac{1+t^k}{1-t^k}.$$

Obtaining a closed form for $G$ is much harder.
The number of circuits of length $n$ is
$$g_n = \sum_{i\in\Z:\,i\equiv 0[k],\,i\equiv n[2]} \binom{n}{\frac{n+i}2},$$
from which, by~\cite[1.54]{gould:identities}, it follows that
$$G(t) = \frac1k\sum_{\zeta^k=1}\frac1{1-\left(\zeta+\zeta^{-1}\right)t}
      = \frac1k\sum_{j=0}^{k-1}\frac1{1-2\cos\left(\frac{2\pi j}k\right)t}.$$

It is not at all obvious how to simplify the above expression. A closed-form
answer can be obtained from~(\ref{eq:c}), namely
$$G(t) = \frac{(2t)^2+\left(1-\sqrt{1-4t^2}\right)^2}{(2t)^2-\left(1-\sqrt{1-4t^2}\right)^2}\cdot\frac{(2t)^k+\left(1-\sqrt{1-4t^2}\right)^k}{(2t)^k-\left(1-\sqrt{1-4t^2}\right)^k},$$
or, expanding,
$$G(t) = \frac{\displaystyle(2t)^k+\sum_{m=0}^{k/2}\left(1-4t^2\right)^m\binom{k}{2m}}{\displaystyle\sum_{m=1}^{(k+1)/2}\left(1-4t^2\right)^m\binom{k}{2m-1}}.$$
However in general this fraction is not reduced. To obtain reduced
fractions for $F(u,t)$ (and thus for $G(t)$), we have to consider
separately the cases where $k$ is odd or even.

For odd $k$, letting $k = 2\ell+1$, we obtain
\begin{gather*}
F(u,t) = \frac{1+(1-u)t}{1-(1+u)t}\cdot\frac{\displaystyle\sum_{m=0}^\ell\alpha_m^\ell(-t)^m(1+(1-u^2)t^2)^{\ell-m}}{\displaystyle\sum_{m=0}^\ell\alpha_m^\ell t^m(1+(1-u^2)t^2)^{\ell-m}},\\
G(t) = \frac{\displaystyle\sum_{m=0}^\ell\alpha_m^\ell(-t)^m}{(1-2t)\left(\displaystyle\sum_{m=0}^\ell\alpha_m^\ell t^m\right)},
\end{gather*}
where
$$\alpha_m^\ell = \begin{cases}\displaystyle(-)^{\frac m2}\binom{\ell-\frac m2}{\frac m2} & \text{if }m\equiv0[2],\\
  \displaystyle(-)^{\frac{m-1}2}\binom{\ell-\frac{m+1}2}{\frac{m-1}2} & \text{if }m\equiv1[2]. \end{cases}$$
For even $k$, with $k = 2\ell$,
\begin{gather*}
F(u,t) = \frac{\displaystyle\sum_{m=0}^{\ell/2}\frac \ell{\ell-m}\binom{\ell-m}{m}(-t^2)^m(1-(1-u^2)t^2)^{\ell-2m}}{(1-(1+u)^2t^2)\left(\displaystyle\sum_{m=0}^{(\ell-1)/2}\binom{\ell-1-m}{m}(-t^2)^m(1-(1-u^2)t^2)^{\ell-1-2m}\right)},\\
G(t) = \frac{\displaystyle\sum_{m=0}^{\ell/2}\frac \ell{\ell-m}\binom{\ell-m}{m}(-t^2)^m}{(1-4t^2)\left(\displaystyle\sum_{m=0}^{(\ell-1)/2}\binom{\ell-1-m}{m}(-t^2)^m\right)},
\end{gather*}
expressed as reduced fractions.

The first few values of $F$, where $\square$ stands for $1+(1-u^2)t^2$, are:
\def\STRUT{\rule[-1ex]{0pt}{3.5ex}}
\begin{center}
\begin{tabular}{r|c@{\quad}r|c}
  $k$ & $F(u,t)$ & $k$ & $F(u,t)$ \\
\hline
\STRUT$1$ & $\frac{1+(1-u)t}{1-(1+u)t}$ &
  $2$ & $\frac{\square}{1-(1+u)^2t^2}$\\[1ex]
$3$ & $\frac{(1+(1-u)t)(\square-t)}{(1-(1+u)t)(\square+t)}$ &
  $4$ & $\frac{\square^2-2t^2}{1-(1+u)^2t^2}$\\[1ex]
$5$ & $\frac{(1+(1-u)t)(\square^2-\square t-t^2)}{(1-(1+u)t)(\square^2+\square t+t^2)}$ &
  $6$ & $\frac{\square^2-3t^2}{(1-(1+u)^2t^2)(\square^2-t^2)}$\\[1ex]
$7$ & $\frac{(1+(1-u)t)(\square^3-\square^2t-2\square t^2+t^3)}{(1-(1+u)t)(\square^3+\square^2t-2\square t^2-t^3)}$ &
  $8$ & $\frac{\square^4-4\square^2t^2+2t^4}{(1-(1+u)^2t^2)(\square^2-2t^2)}$\\[1ex]
$9$ & $\frac{(1+(1-u)t)(\square-t)(\square^3-3\square t^2-t^3)}{(1-(1+u)t)(\square+t)(\square^3-3\square t^2+t^3)}$ &
  $10$ & $\frac{\square^4-5\square^2t^2+5t^4}{(1-(1+u)^2t^2)(\square^4-3\square^2t^2+t^4)}$\\[1ex]
$11$ & $\frac{(1+(1-u)t)(\square^5-\square^4t-4\square^3t^2+3\square^2t^3+3\square t^4-t^5)}{(1-(1+u)t)(\square^5+\square^4t-4\square^3t^2-3\square^2t^3+3\square t^4+t^5)}$ &
  $12$ & $\frac{(\square^2-2t^2)(\square^4-4\square^2t^2+t^4)}{(1-(1+u)^2t^2)(\square^2-3t^2)(\square^2-t^2)}$
\end{tabular}
\end{center}
These rational expressions were computed and simplified using the
computer algebra program \textsc{Maple}.

\subsection{Trees}\label{subs:trees}
Let $\gf$ be the $d$-regular tree. Then
$$F(0,t)=1$$ as a tree has no proper circuit; while direct (i.e.,
without using Corollary~\ref{cor:main}) computation of $G$ is more
complicated.  It was first performed by Kesten~\cite{kesten:rwalks};
here we will derive the extended circuit series $F(u,t)$ and also
obtain the answer using Corollary~\ref{cor:main}.

Let $\tree$ be a regular tree of degree $d$ with a fixed root $\birth$, and
let $\tree'$ be the connected component of $\birth$ in the two-tree forest
obtained by deleting in $\tree$ an edge at $\birth$. Let $F(u,t)$ and
$F'(u,t)$ respectively count circuits at $\birth$ in $\tree$ and $\tree'$.
For instance if $d=2$ then $F'$ counts circuits in $\N$ and $F$ counts
circuits in $\Z$.  For a reason that will become clear below, we make
the convention that the empty circuit is counted as `$1$' in $F$ and
as `$u$' in $F'$. Then we have
\begin{gather*}
F' = u + (d-1)tF't\frac1{1-(d-2+u)tF't},\\
F = 1 + dtF't\frac1{1-(d-1+u)tF't}.\\
\end{gather*}
Indeed a circuit in $\tree'$ is either the empty circuit (counted as
$u$), or a sequence of circuits composed of, first, a step in any of
$d-1$ directions, then a `subcircuit' not returning to $\birth$, then
a step back to $\birth$, followed by a step in any of $d-1$ directions
(counting an extra factor of $u$ if it was the same as before), a
subcircuit, etc. If the `subcircuit' is the empty circuit, it
contributes a bump, hence the convention on $F'$.  Likewise, a circuit
in $\tree$ is either the empty circuit (now counted as $1$) or a
sequence of circuits in subtrees each isomorphic to $\tree'$.

We solve these equations to
\begin{gather*}
F'(1-u,t) = \frac{2(1-u)}{1-u(d-u)t^2+\sqrt{(1 + u(d-u)t^2)^2-4(d-1)t^2}},\\
F(1-u,t) = \frac{2(d-1)(1-u^2t^2)}{(d-2)(1 + u(d-u)t^2)+d\sqrt{(1 + u(d-u)t^2)^2-4(d-1)t^2}}.
\end{gather*}

Using~(\ref{eq:c}) and $F(0,t)=1$ we would obtain
$$G(t) = \frac{1+(d-1)\left(\frac{1-\sqrt{1-4(d-1)t^2}}{2(d-1)t}\right)^2}{1-\left(\frac{1-\sqrt{1-4(d-1)t^2}}{2(d-1)t}\right)^2},$$
or, after simplification,
$$G(t) = \frac{2(d-1)}{d-2+d\sqrt{1-4(d-1)t^2}},$$
which could have been obtained by setting $u=0$ in $F(1-u,t)$.

In particular if $d=2$, then $\gf=C_\infty=\Z$ and
$$G(t) = \sum_{n\ge0}\binom{2n}{n}t^{2n} = \frac1{\sqrt{1-4t^2}}.$$

Note that for all $d$ the $d$-regular tree $\gf$ is the Cayley graph
of $\Gamma=(\Z/2\Z)^{*d}$ with standard generating set. If $d$ is even,
$\gf$ is also the Cayley graph of a free group of rank $d/2$ generated
by a free set.  We have thus computed the spectral radius of a random
walk on a freely generated free group: it is, for $(\Z/2\Z)^{*d}$ or for
$\free{d/2}$, equal to
\begin{equation}\label{specrad:tree}
  \frac{2\sqrt{d-1}}{d}.
\end{equation}

Note that for $d=2$ the series $F(u,t)$ does have a simple
expansion. By direct expansion, we obtain the number of circuits of
length $2n$ in $\Z$, with $m$ local extrema, as
$$(t^{2n}u^m|F(u,t)) =
  \begin{cases} 2\binom{n-1}{\frac{m-1}2}^2 & \text{if }m\equiv1[2], \\
                2\binom{n-1}{\frac{m}{2}}\binom{n-1}{\frac{m-2}{2}} & \text{if }m\equiv0[2].
  \end{cases}$$

We may even look for a richer generating series than $F$: let
$$H(u,v,t) = \sum_{\pi:\text{ path starting at }\birth}u^{\bcount(\pi)}v^{\dist(\birth,\pi_{|\pi|})}t^{|\pi|}\in\N[u,v][[t]],$$
where $\dist$ denotes the graph distance. Then
\begin{align*}
H(1,v,t) &= F(1,t) + dF'tvF + dF'tv(d-1)F'tvF + \dots\\
  &= \frac{1+F'(1,t)tv}{1-(d-1)F'(1,t)tv}F(1,t);
\end{align*}
and as $H$ is a sum of series counting paths between fixed vertices
we obtain $H(u,v,t)$ from $H(1,v,t)$ by extending~(\ref{eq:b}) linearly:
$$\frac{H(1-u,v,t)}{1-u^2t^2} = \frac{H\left(1,v,\frac{t}{1+u(d-u)t^2}\right)}{1+u(d-u)t^2}.$$
We could also have started by computing
$$H(0,v,t) = \frac{1+vt}{1-(d-1)vt},$$
the growth series of all proper paths in $\tree$,
and using~(\ref{eq:c}) and~(\ref{eq:d}) obtain
\begin{gather*}
  \def\xxx{1-\sqrt{1-4(d-1)t^2}}
  H(1,v,t) = \frac{1+\left(\frac\xxx{2t}\right)^2}{1-u^2\left(\frac\xxx{2(d-1)t}\right)^2}\cdot H\left(\frac\xxx{2(d-1)t},0,v\right),\\
  H(u,v,t) = \frac{1-t^2u^2}{1+u(d-u)t^2}\cdot\frac{(d-1)(4t^2+\square^2)}{4(d-1)^2t^2-u^2\square^2}\cdot\frac{2(d-1)t+v\square}{2t-v\square},\\
\end{gather*}
where $\square = 1+u(d-u)t^2-\sqrt{(1+u(d-u)t^2)^2-4(d-1)t^2}$.

Recall that the growth series of a graph $\gf$ at a base point $\birth$ is
the power series
$$P(t) = \sum_{v\in V(\gf)}t^{\dist(\birth,v)},$$
where $\dist$ denotes the distance in $\gf$. The series $H$ is very
general in that it contains a lot of information on $\tree$, namely
\begin{itemize}
\item $H(u,0,t)=F(u,t)$;
\item $H(0,1,t)= \frac{1+t}{1-(d-1)t} = P(t)$ is the growth series of $\tree$;
\item $H(1,1,t) = 1/(1-dt)$ is the growth series of all paths in $\tree$.
\end{itemize}
(Note that these substitutions yield well-defined series because for
any $i$ there are only finitely many monomials having $t$-degree equal
to $i$.)

We can also use this series $H$ to compute the circuit series
$F_C$ of the cycle of length $k$, that was found in the
previous section. Indeed the universal cover of a cycle is the regular
tree $\tree$ of degree $2$, and circuits in $C$ correspond bijectively to
paths in $\tree$ from $\birth$ to any vertex at distance a multiple of
$k$. We thus have
$$F_C(u,t) = \sum_{\zeta: \zeta^k=1} H(u,\zeta,t)$$
where the sum runs over all $k$th roots of unity and $d=2$ in $H$.

We consider next the following graphs: take a $d$-regular tree and fix
a vertex $\birth$. At $\birth$, delete $e$ vertices and replace them by $e$
loops. Then clearly
$$F(0,t) = \frac{1+t}{1-(e-1)t},$$
as all the non-backtracking paths are constrained to the $e$
loops. Using~(\ref{eq:c}), we obtain after simplifications
\begin{equation}\label{eq:treeloop}
  G(t) = \frac{2(d-1)}{d+e-2-2e(d-1)t + (d-e)\sqrt{1-4(d-1)t^2}}.
\end{equation}
The radius of convergence of $G$ is
$$\min\left\{\frac1{2\sqrt{d-1}},\frac{e-1}{d+e^2-2e}\right\}.$$

\subsection{Tougher Examples}\label{ex:ladder}
\let\L=\leftrightarrow \let\U=\uparrow \let\D=\downarrow In this
subsection we outline the computations of $F$ and $G$ for more
complicated graphs. They are only provided as examples and are
logically independent from the remainder of the paper. The arguments
will therefore be somewhat condensed.

First take for $\gf$ the Cayley graph of $\Gamma=(\Z/2\Z)\times\Z$
with generators $(0,-1)=\text{`$\D$'}$, $(0,1)=\text{`$\U$'}$ and
$(1,0)=\text{`$\L$'}$. Geometrically, $\gf$ is a doubly-infinite
two-poled ladder.

In Subsection~\ref{subs:trees} we computed
$$F_\Z(u,t) = \frac{1-(1-u)^2t^2}{\sqrt{(1+(1-u^2)t^2)^2-4t^2}},$$
the growth of circuits restricted to one pole of the ladder.  A
circuit in $\gf$ is a circuit in $\Z$, before and after each step
($\U$ or $\D$) of which we may switch to the other pole (with a $\L$)
as many times as we wish, subject to the condition that the circuit
finish at the same pole as it started. This last condition is expressed
by the fact that the series we obtain must have only coefficients of
even degree in $t$. Thus, letting $\even(f) = \frac{f(t)+f(-t)}{2}$,
we have
$$G(t) = \even\left(\frac1{1-t} F_\Z\left(1,\frac{t}{1-t}\right)\right);$$
it is then simple to obtain $F(u,t)$ by performing the
substitution~(\ref{eq:c}).

The following direct argument also gives $F(u,t)$: a walk on the
ladder is obtained from a walk on a pole (i.e.\ on $\Z$) by inserting
before and after every step on a pole a (possibly empty) sequence of
steps from one pole to the other. This process is expressed by
performing on $F_\Z$ the substitution
$$t\mapsto t + t^2 + t^3u + t^4u^2 + \dots = t + \frac{t^2}{1-tu},$$
corresponding to replacing a step on a pole by itself, or itself
followed by a step to the other pole, or itself, a step to the other
pole and a step back, etc. But if the path had a bump at the place the
substitution was performed, this bump would disappear when a step is
added from one pole to the other. In formulas,
$$tu\mapsto tu + t^2 + t^3u + t^4u^2 + \dots = tu + \frac{t^2}{1-tu}.$$
Finally we must add at the beginning of the path a sequence of steps
from one pole to the other. Therefore we obtain
$$F(u,t) = \even\left\{\left(1+\frac{t}{1-tu}\right)F_\Z\left(\frac{tu+t^2/(1-tu)}{t+t^2/(1-tu)},t+\frac{t^2}{1-tu}\right)\right\}.$$

As another example, consider the group $\Z$ generated by the non-free
set $\{\pm1,\pm2\}$. Geometrically, it can be seen as the set of
points $(2i,0)$ and $(2i+1,\sqrt{3})$ for all $i\in\Z$, with edges
between all points at Euclidean distance $2$ apart; but we will not
make use of this description. The circuit series of $\Z$ with this
enlarged generating set will be an algebraic function of degree $4$
over the rationals.

Define first the following series:
\begin{description}
\item[$f(t)$] counts the walks from $0$ to $0$ in $\N$;
\item[$g(t)$] counts the walks from $0$ to $1$ in $\N$;
\item[$h(t)$] counts the walks from $1$ to $1$ in $\N$.
\end{description}
\let\U=\uparrow
\let\D=\downarrow
\let\UU=\upuparrows
\let\DD=\downdownarrows
Denote the generators of $\Z$ by $1=\U$, $2=\UU$, $-1=\D$ and~$-2=\DD$. 
The series then satisfy the following equations, where the generators'
symbol is written instead of `$t$' to make the formulas self-explanatory:
\begin{gather*}
  f = 1 + \left(\U f\D + \U g\DD + \UU g\D + \UU h\DD\right)f,\\
  g = f\U f + f\UU g,\\
  h = f + f\D g + g\DD g,
\end{gather*}
giving a solution $f$ that is algebraic of degree $4$ over $\Z(t)$.

Then define the following series:
\begin{description}
\item[$G$] counts the walks from $0$ to $0$ in $\Z$;
\item[$e$] counts the walks from $0$ to $1$ in $\Z$.
\end{description}
They satisfy the equations
\begin{gather*}
  G = 1 + 2\left(\U f\D G + \UU g\D G + \U f\DD e + \UU g\DD e + \U g\DD G + \UU h\DD G\right),\\
  e = G\U f + G\UU f + G\UU g
\end{gather*}
giving the solution

$$G = \frac{4 + 3t - 6t^2 - 10t(1+2t)\delta + 2t^2(3+8t)\delta^2 - 6t^4(1+t)\delta^3}{4 - 7t - 36t^2}$$
where $\delta$ is a root of the equation
$$1 - (2t+1)\delta + t(2+3t)\delta^2 - t^2(1+2t)\delta^3 + t^4\delta^4 = 0.$$

\section{Cogrowth of Non-Free Presentations}\label{sec:nonfree}
We perform here a computation extending the results of
Section~\ref{subs:gptheory}. The general setting, expressed in the
language of group theory, is the following: let $\Pi$ be a group
generated by a finite set $S$ and let $\Xi<\Pi$ be any subgroup.
We consider the following generating series:
\begin{align*}
  F(t) &= \sum_{\gamma\in\Xi<\Pi}t^{|\gamma|},\\
  G(t) &= \sum_{\substack{\text{words $w$ in $S$}\\ \text{defining an
        element in $\Xi$}}}t^{|w|},
\end{align*}
where $|\gamma|$ is the minimal length of $\gamma$ in the generators
$S$, and $|w|$ is the usual length of the word $w$. Is there some
relation between these series? In case $\Pi$ is quasi-free on $S$, the
relation between $F$ and $G$ is given by Corollary~\ref{cor:main}.  We
consider two other examples: $\Pi$ quasi-free but on a set smaller
than $S$, and $\Pi = PSL_2(\Z)$.

\subsection{$\boldsymbol\Pi$ Quasi-Free}
Let $S$, $T$ be finite sets, and $\overline\cdot$ an involution on
$S$. Consider the two presentations
\begin{gather*}
  \Pi = \langle S|\,s\overline s=1\,\forall s\in S\rangle,\\
  \Pi = \langle S\cup T|\,s\overline s=1\,\forall s\in S;\;t=1\,\forall t\in T\rangle.
\end{gather*}
Let $\Xi<\Pi$ be any subgroup, and let $F'$ and~$G'$ be
the generating series related to the first presentation. Clearly
$F'=F$, as both series count the same objects in $\Pi$
(regardless of $\Pi$'s presentation); while
$$G(t) = \frac{G'\left(\frac{t}{1-|T|t}\right)}{1-|T|t}.$$
Indeed any word $w=w_1\dots w_n$ in $S\cup T$ defining an element of
$\Xi$ can be uniquely decomposed as $w=t_0\,s_1t_1\dots s_mt_m,$
where $s_i\in S$, $t_i$ are words in $T$ for all $i$, and $s_1\dots
s_n$ defines an element of $\Xi$; moreover all choices of $s_1\dots
s_n$ defining an element of $\Xi$ and words $t_i$ in $T$ give a
distinct word $w$. It then suffices to note that the generating series
for any of the $t_i$ is $1/(1-|T|t)$.

Putting everything together, we obtain:
\begin{prop}
  Let $\Pi$ be as above, $\Xi<\Pi$ a subgroup.
  Then
  $$\frac{F(t)}{1-t^2} = \frac{G\left(\frac{t}{1+|T|t+(|S|-1)t^2}\right)}{1+|T|t+(|S|-1)t^2}.$$
\end{prop}

\subsection{$\boldsymbol{\Pi = PSL_2(\Z)}$}
Let $$\Pi = PSL_2(\Z) = \langle a,b|\,a^2,b^3\rangle,$$
and let $\Xi<\Pi$ be any subgroup. We take $S=\{a,b,b^{-1}\}$.

We suppose $\Xi$ is torsion-free, i.e.\ contains no element of the
form $waw^{-1}$ or $wb^{\pm1}w^{-1}$. Let $\gf$ be the Schreier graph
of $(\Pi,\{a,b,b^{-1}\})$ relative to $\Xi$, as defined in
Subsection~\ref{subs:gptheory}; this is a trivalent graph whose vertex
set is $\Xi\backslash\Pi$. Its vertices can be grouped in triples
$w^\Delta = \{w,wb,wb^{-1}\}$ connected in triangles. Let $\gff$ be
the graph obtained from $\gf$ by identifying each triple to a vertex.
Explicitly,
\begin{gather*}
  V(\gff) = \{w^\Delta:\,w\in V(\gf)\},\\
  E(\gff) = \{(v^\Delta,(va)^\Delta): v\in V(\gf)\};
\end{gather*}
the involution on $E(\gff)$ is the switch $(A,B)\mapsto(B,A)$ and the
extremity functions $E(\gff)\to V(\gff)$ are the natural projections.
Note that $\gff$ is a $3$-regular graph (for instance, $1^\Delta$ is
connected to $a^\Delta$, $(ba)^\Delta$ and $(b^{-1}a)^\Delta$). In case
$\Xi=1$, it is the $3$-regular tree. By construction we have a
$3$-to-$1$ map $\Delta:V(\gf)\to V(\gff)$.  We fix an origin $\birth=1^\Delta$
in $\gff$, and let $F_\gff(u,t)$ be the circuit series of $(\gff,\birth)$.

Let $\gfe$ be a triangle, $G_\gfe(t)$ count the circuits at a fixed
vertex of $\gfe$ and $G_\gfe^{\neq}(t)$ count paths between two fixed
distinct vertices of $\gfe$. These series were computed in
Section~\ref{subs:complete}, with $G_\gfe^{\neq}(t) = F'(1,t) + F''(1,t)$.

Circuits at $\birth$ in $\gf$ can be projected to circuits at $\birth$ in $\gff$
simply by deleting all edges of type $(w,wb^{\pm1})$ and projecting
the other edges through $\Delta$. Conversely, circuits in $\gff$ can
be lifted to $\gf$ by lifting the edges through $\Delta^{-1}$, and
connecting them in $\gf$ with arbitrary paths remaining inside the
triples; to lift the path $\pi=(\pi_1,\dots,\pi_n)$ from $\gff$ to
$\gf$, choose edges $\rho_1,\dots,\rho_n$ with
$(\rho_i^\alpha)^\Delta=\pi_i^\alpha$ and
$(\rho_i^\omega)^\Delta=\pi_i^\omega$ for all $i\in\{1,\dots,n\}$, and
choose, for all $i\in\{0,\dots,n\}$, paths $\tau_i$ from
$\rho_i^\omega$ to $\rho_{i+1}^\alpha$ remaining inside
$(\rho_i^\omega)^\Delta$, where by convention
$\rho_0^\omega=\rho_{n+1}^\alpha=\birth$. Then the lift corresponding to
these choices is
\begin{equation}\label{eq:shuffle}
  \tau_0\cdot\rho_1\cdot\tau_1\cdots\rho_n\cdot\tau_n.
\end{equation}
Furthermore all circuits at $\birth$ in $\gf$ can be obtained this way.

Define $\overline G$ as the series counting paths that start and
finish at a vertex in the same triple as $\birth$.
It can be obtained using~(\ref{eq:shuffle}) by letting $\rho$ range
over all paths in $\gff$, and for each choice of $\rho$ and for each
$i\in\{1,\dots,n-1\}$ letting $\tau_i$ range over $G_\gfe$ or
$G_\gfe^{\neq}$ depending on whether $\rho$ has or not a bump at
$i$, and letting $\tau_0$ and $\tau_n$ range over all paths inside the
triple $\birth^\Delta$. In equations, this relation is expressed as
$$\overline G(t) = \left(\frac1{1-2t}\right)^2/G_\gfe(t)\cdot F_\gff\left(G_\gfe^{\neq}(t)/G_\gfe(t),tG_\gfe(t)\right).$$
Now the series $G$ we wish to obtain is approximately $\overline
G(t)/9$: for any choice of $x,y\in \birth^\Delta$ there are
approximately the same number of long enough paths from $x$ to $y$.

A summand of $F(t)$ is the unique lifting of a summand of
$F_\gff(0,t)$, but is twice longer in $\gf$ than in $\gff$.

\begin{defn}
  Two series $A(t)$, $B(t)$ are \emdef{equivalent}, written $A\sim B$,
  if they have the same radius of convergence $\rho$, and
  there exists a constant $K$ such that
  $$\frac1K<A(t)/B(t)<K\text{ as }t\to\rho.$$
\end{defn}
Then the remarks of the previous paragraph can be written as
\begin{gather*}
  F(t) \sim F_\gff(0,t^2),\\
  G(t) \sim F_\gff(G_\gfe^{\neq}(t)/G_\gfe(t),tG_\gfe(t)).
\end{gather*}
Letting $G_\gff$ be the circuit series of $\gff$, we use
Corollary~\ref{cor:main} to obtain
\begin{gather*}
  G_\gfe(t)^{\neq} = \frac{t}{1-t-2t^2},\qquad G_\gfe(t) = \frac{1-t}{1-t-2t^2},\\
  F(t) \sim G_\gff\left(\frac{t^2}{1+2t^4}\right),\\
  G(t) \sim G_\gff\left(\frac{t^2}{1-t-3t^2}\right)\\
\intertext{so}
  F(t) \sim G\left(\frac{t\sqrt{4+13t^2-8t^4}-t^2}{2(1+t^2)(1+2t^2)}\right).
\end{gather*}

\begin{figure}
\begin{center}
\setlength\unitlength{1pt}
\begin{picture}(300,300)
\put(0,170){\makebox(0,0)[lb]{\smash{\normalsize $\nu$}}}
\put(150,30){\makebox(0,0)[lb]{\smash{\normalsize $\alpha$}}}
\put(-30,-10){\epsfig{file=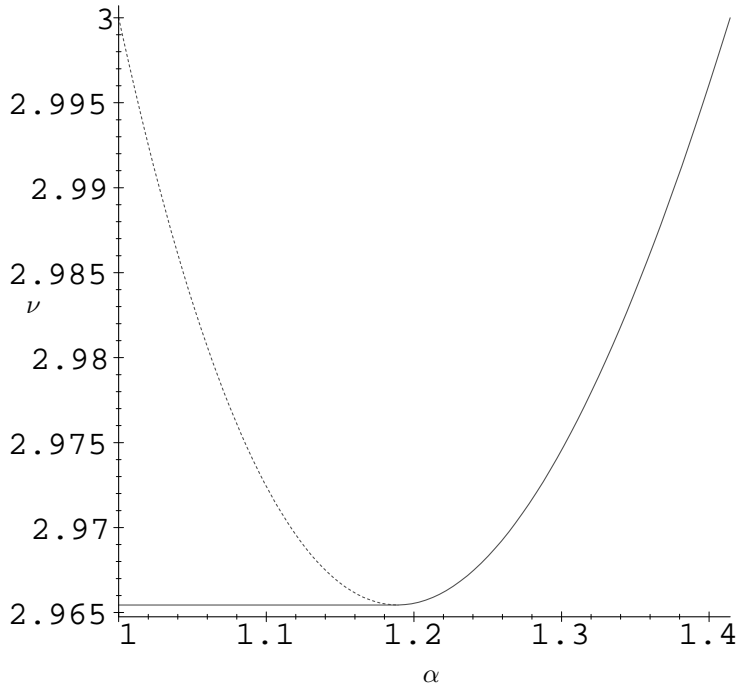}}
\end{picture}
\end{center}
\caption{The function $\alpha\mapsto\nu$ relating cogrowth and spectral
  radius, for subgroups of $PSL_2(\Z)$.}
\end{figure}

Let $\gf$ be a simplicial complex such that at each vertex an edge and
a (filled-in) triangle meet; choose a base point $\birth$ in $\gf$.  Say a
circuit in the $1$-skeleton of $\gf$ is \emdef{reduced} if it contains
no bump nor two successive edges in the same triangle; thus reduced
circuits are in bijection with homotopy classes of $(\gf,\birth)$. Let
$F(t)$ be the proper circuit series and $G(t)$ the circuit series of
$\gf$. Let
$$(t)\phi = \frac{t\sqrt{4+13t^2+8t^4}-t^2}{2(1+t^2)(1+2t^2)}.$$

We have proved the following theorem and corollary, similar to those
in Section~\ref{subs:gptheory}:
\begin{thm} $F(t)\sim G((t)\phi)$.
\end{thm}

\begin{cor}
  Let $\Xi$ be a subgroup of $\Pi=PSL_2(\Z)$; let $\nu$ be the
  spectral radius of the simple random walk on
  $\Xi\backslash\Pi$, and $\alpha$ the ``cogrowth'' rate of
  $\Xi\backslash\Pi$.  Then provided that
  $\alpha\in[\sqrt\rho,\rho]$, where $\rho$ is the word growth of
  $\Pi$, namely $\sqrt2$, we have
  $$1/\nu = (1/\alpha)\phi,$$
  so
  $$\nu = \frac12\sqrt{8\alpha^{-2}+13+4\alpha^2}+\frac12.$$
\end{cor}

\begin{proof}
  The function $\phi$ is monotonously increasing between $0$ and
  $1/\sqrt[4]{2}$, where it reaches its maximum. The same argument
  applies as that given in the proof of Corollary~\ref{cor:homo}.
\end{proof}

We now state the same results for an arbitrary virtually free group
with an appropriate generating system. Let $\Pi$ be a virtually
free group, such that there is a split exact sequence
$$\begin{diagram}
  \node{1} \arrow{e} \node{\Sigma} \arrow{e} \node{\Pi} \arrow{e,t}{\pi} \node{\Upsilon} \arrow{e} \arrowm{0}{-150}\arrow{w,m} \node{1}
\end{diagram}$$
where $\Upsilon$ is a finite group and $\Sigma$ has a presentation
$$\Sigma = \langle s\in S|\,s\overline s=1\;\forall s\in S\rangle.$$
We assume further that $\Pi$ is generated by a set $T=T'\sqcup T''$
with $T''$ in bijection through $\pi$ with $\Upsilon\setminus\{1\}$,
$T'$ mapping through $\pi$ to $\{1\}$, and $T'\times(T''\cup\{1\})$ in
bijection with $S$ through $(t,p)\mapsto p^{-1}tp$.

For example, consider $\Pi=PSL_2(\Z)=\langle a,b,b^{-1}\rangle$. Take
$T'=\{a\}$ and $T''=\{b,b^{-1}\}$, take $\Upsilon=\langle
b,b^{-1}\rangle$ and $\Sigma=\langle a,bab^{-1},b^{-1}ab\rangle$. Then
the hypotheses are satisfied.

With these hypotheses, the Cayley graph $\gf$ of $\Pi$ is a
collection of complete graphs of size $|\Upsilon|$, with at each vertex
$|T'|$ edges leaving to other complete graphs, and such that if each
of these complete graphs is shrunk to a point the resulting graph is a
tree. The following theorem is then a straightforward generalization
of the argument given for $PSL_2(\Z)$.

\begin{thm}
  With the notation introduced above, let $\Xi$ be any subgroup of
  $\Pi$ not intersecting $\{t^\gamma|\,t\in T,\gamma\in\Pi$ and let
  $F(t)$, $G(t)$ be the ``cogrowth'' series and circuit series of
  $\Xi\backslash\Pi$.  Let $\gfe$ be the complete graph on
  $|\Upsilon|$ vertices and let $G_\gfe(t)$, $G_\gfe^{\neq}(t)$ count
  the circuits and the non-closing paths respectively in $\gfe$.
  Define the function $\phi$ by
  $$\left(\frac{t^2}{1+(|S|-1)t^4}\right)\phi = \frac{tG_\gfe}{1+(G_\gfe-G_\gfe^{\neq})((|S|-1)G_\gfe+G_\gfe^{\neq})t^2}.$$

  Then we have
  $$F(t) \sim G((t)\phi).$$
\end{thm}

\section{Free Products of Graphs}\label{sec:freeamalg}
We give here a general construction combining two pointed graphs and
show how to compute the generating functions for circuits in the
``product'' in terms of the generating functions for circuits in the
factors.

\def\0{\cite[Definition~4.8]{quenell:freeprod}}
\begin{defn}[Free Product, \0] Let
  $(\gfe,\birth)$ and $(\gff,\birth)$ be two connected pointed graphs. Their
  \emdef{free product} $\gfe*\gff$ is the graph constructed as
  follows: start with copies of $\gfe$ and $\gff$ identified at $\birth$;
  at each vertex $v$ apart from $\birth$ in $\gfe$, respectively $\gff$,
  glue a copy of $\gff$, respectively $\gfe$, by identifying $v$ and
  the $\birth$ of the copy. Repeat the process, each time gluing $\gfe$'s
  and $\gff$'s to the new vertices.

\end{defn}

If $(E,S)$, $(F,T)$ are two groups with fixed generators whose Cayley
graphs are $\gfe$ and $\gff$ respectively, then $\gfe*\gff$ is the
Cayley graph of $(E*F,S\sqcup T)$. 

%
%

We will now compute the circuit series of $\gfe*\gff$ in terms of the
circuit series of $\gfe$ and~$\gff$.  Let $G_\gfe$, $G_\gff$
and~$G_\gf$ be the generating functions counting circuits in $\gfe$,
$\gff$ and $\gf=\gfe*\gff$ respectively. We will use the following
description: given a circuit at $\birth$ in $\gf$, it can be decomposed as
a product of circuits never passing through $\birth$. Each of these
circuits, in turn, starts either in the $\gfe$ or the $\gff$ copy at
$\birth$. Say one starts in $\gfe$; it can then be expressed as a circuit in
$\gfe$ never passing through $\birth$, and such that at all vertices,
except the first and last, a circuit starting in $\gff$ has been
inserted. Moreover, any choice of such circuits satisfying these
conditions will give a circuit at $\birth$ in $\gf$, and different
choices will yield different circuits.

Let $H_\gfe$ (respectively $H_\gff$) be the generating
function counting non-trivial circuits in $\gfe$ (respectively $\gff$)
never passing through $\birth$.  Obviously
$$G_\gfe = \frac1{1-H_\gfe}\qquad\text{so }H_\gfe = 1 - \frac1{G_\gfe}.$$
Let $L_\gfe$ (respectively $L_\gff$) be the generating
function counting non-trivial circuits in $\gf$ never passing through
$\birth$ and starting in $\gfe$ (respectively $\gff$). Then
\begin{align*}
  L_\gfe(t) &= H_\gfe\left(\frac{t}{1-L_\gff(t)}\right)\cdot(1-L_\gff(t)),\\
  L_\gff(t) &= H_\gff\left(\frac{t}{1-L_\gfe(t)}\right)\cdot(1-L_\gfe(t)).
\end{align*}
Indeed write $H_\gfe = \sum h_n t^n$. Then by the description
given above
$$L_\gfe = \sum h_n t^n \left(\frac1{1-L_\gff}\right)^{n-1},$$
which is precisely the given formula. Finally
$$G_\gf = \frac1{1-L_\gfe-L_\gff}.$$

Writing $M_\gfe=t/(1-L_\gfe)$, and similarly for $\gff$, we
simplify these equations to
\begin{gather}
  1-\frac{t}{M_\gfe} = \left(1-\frac1{G_\gfe(M_\gff)}\right)\cdot\frac{t}{M_\gff},\notag\\
  1-\frac{t}{M_\gff} = \left(1-\frac1{G_\gff(M_\gfe)}\right)\cdot\frac{t}{M_\gfe},\notag\\
  G_\gf = \frac1{1 - \left(1 - \frac{t}{M_\gfe}\right) - \left(1 - \frac{t}{M_\gff}\right)} = \frac{1/t}{1/M_\gfe + 1/M_\gff - 1/t},\notag\\
\intertext{so}
  \frac1{M_\gfe} + \frac1{M_\gff} - \frac1t = \frac1{M_\gfe G_\gff(M_\gfe)} = \frac1{M_\gff G_\gfe(M_\gff)} = \frac1{tG_\gf}.\label{eq:free1}
\end{gather}
If $f$ is a power series with $f(0)=0$ and $f'(0)\neq0$, let us write
$f^{-1}$ for the inverse of $f$; i.e.\ for the series $g$ with
$g(f(t))=f(g(t))=t$ (for instance, $f(t)=t$ is equal to its inverse;
the inverse of $\frac{at+b}{ct+d}$ is $\frac{dt-b}{-ct+a}$).

From~(\ref{eq:free1}) we obtain $M_\gfe = (tG_\gff)^{-1}\circ(tG_\gf)$
and $M_\gff = (tG_\gfe)^{-1}\circ(tG_\gf)$; so composing
(\ref{eq:free1}) with $(tG_\gf)^{-1}$ we obtain the
\begin{thm}\label{thm:free}
  \begin{equation}\label{eq:free}
    \frac1{(tG_\gf)^{-1}} = \frac1{(tG_\gfe)^{-1}} + \frac1{(tG_\gff)^{-1}} - \frac1t.
  \end{equation}
\end{thm}
An equation equivalent to this one, though not trivially so, appeared
in a paper by Gregory Quenell~\cite{quenell:freeprod}, and, in yet another
language, in a paper by Dan Voiculescu~\cite[Theorem~4.5]{voiculescu:rv}.

We can use~(\ref{eq:free}) to obtain by a different method
the circuit series of regular trees (see Section~\ref{subs:trees}).
Indeed the free product of regular trees of degree $d$ and $e$ is a
regular tree of degree $d+e$. Letting $G_d$ denote the circuit series
of a regular tree of degree $d$, we ``guess'' that
$$G_d(t) = \frac{2(d-1)}{d-2+d\sqrt{1-4(d-1)t^2}},$$
and verify that the limit
$$G_1(t) = \lim_{d\to1}\frac{2(d-1)}{d-2+d\sqrt{1-4(d-1)t^2}}=\frac1{1-t^2}$$
is indeed the circuit series of the $1$-regular tree.
Then we compute
$$(tG_d)^{-1}(u) = \frac{2u}{2-d+d\sqrt{1+4u^2}};$$
if we define $\triangle_d$ by
$$\triangle_d := \frac1{(tG_d)^{-1}(u)}-\frac1u =
d\frac{\sqrt{1+4u^2}-1}{2u},$$
it satisfies $\triangle_d+\triangle_e=\triangle_{d+e}$ and we have proved
that our guess of $G_d$ is correct, in light of~(\ref{eq:free}).

As another application of~(\ref{eq:free}), we will compute the circuit
series $G(t)$ of the Cayley graph of $PSL_2(\Z)=\langle
a,b|\,a^2,b^3\rangle$ with generators $\{a,b,b^{-1}\}$. This graph is
the free product of the $1$-regular tree $\gfe$ and of the $3$-cycle $\gff$.
We know from Section~\ref{subs:complete} that
$$G_\gfe=\frac1{1-t^2},\qquad G_\gff=\frac{1-t}{(1+t)(1-2t)}$$
are the circuit series of $\gfe$ and $\gff$.
We then compute
\begin{gather*}
  (tG_\gfe)^{-1}(u) = \frac{\sqrt{1+4u^2}-1}{2u},\\
  (tG_\gff)^{-1}(u) = \frac{1+u-\sqrt{1-2u+9u^2}}{2(1-2u)},
\end{gather*}
so after some lucky simplifications
\begin{align*}
  G(t) &= \frac1t\left(\frac1{1/(tF_\gfe)^{-1}(u)+1/(tG_\gff)^{-1}(u)-1/u}\right)^{-1}\\
  &= \frac{(2-t)\sqrt{1-2t-5t^2+6t^3+t^4}-t+t^2+t^3}{2(1-2t-5t^2+6t^3)}.
\end{align*}
(A closed form such as this one does not exist for $(\Z/2\Z)*(\Z/k\Z)$ with
$k>3$, because then the series $G$ is algebraic of degree greater than
2.)

\begin{cor}
  If the circuit series of $\gfe$ and $\gff$ are both algebraic, then
  the circuit series of $\gfe*\gff$ is also algebraic.
\end{cor}

\begin{proof}
  Sums and products of algebraic series are algebraic. If $f$
  satisfies the algebraic relation $P(f,t)=0$, then its formal inverse
  satisfies the relation $P(t,f^{-1})=0$ so is also algebraic.
\end{proof}

\begin{lem}
  We have
  $$\rho(f) = \sup_{t}(tf(t))^{-1},$$
  where the supremum is taken over all $t$ such that the series
  $(tf)^{-1}$ converges.
  If $f$ is $\rho$-recurrent, then also
  $$\rho(f) = \lim_{t\to\infty}(tf(t))^{-1}.$$
\end{lem}

\begin{proof}
  Clearly $\rho(f)=\rho(tf)$; if $tf$ converges over $[0,\rho[$ then
  $(tf)^{-1}$ converges over $[0,\sigma[$ where $\sigma=\rho
  f(\rho)$; then $\lim_{t\to\sigma}(tf)^{-1}=\rho$. The second
  assertion follows because in this case $\sigma=\infty$.
\end{proof}

\begin{cor}
  Let the circuit series of $\gfe$, $\gff$ and~$\gf=\gfe*\gff$ be
  $G_\gfe$, $G_\gff$ and $G_\gf$ respectively, and suppose all three
  series are recurrent. Then
  $$1/\rho(G_\gf) = 1/\rho(G_\gfe) + 1/\rho(G_\gff).$$
\end{cor}
\begin{proof}
  This follows from
  \begin{align*}
    1/\rho(G_\gf) &= \lim_{t\to\infty}\frac1{(tG_\gf)^{-1}}\\
    &= \lim_{t\to\infty}\frac1{(tG_\gfe)^{-1}}+\frac1{(tG_\gff)^{-1}}-\frac1t\qquad\text{by~(\ref{eq:free})}\\
    &= 1/\rho(G_\gfe) + 1/\rho(G_\gff) - 0.
  \end{align*}
\end{proof}
Note that the corollary does not extend to non-recurrent series; for
instance, it fails if $\gfe=\gff=\Z$. Indeed then
\begin{gather*}
  G_\gfe = G_\gff = \frac1{\sqrt{1-4t^2}},\qquad\rho(G_\gfe)=\rho(G_\gff)=1/4,\\
  F_\gf = \frac3{1+2\sqrt{1-12t^2}},\qquad\rho(G_\gf)=1/\sqrt{12}.
\end{gather*}

\section{Direct Products of Graphs}\label{sec:direct}
There are two natural definitions for \emdef{direct products} of
graphs; they correspond to direct products of groups with generating
set either the union or cartesian product of the generating sets of the
factors. A general treatment of products of graphs can be found
in~\cite[pages 65 and 203]{cvetkovic:spectra}

\begin{defn}
  If $S$ is a set, the \emdef{stationing graph} on $S$ is the graph
  $\gf=\Sigma_S$ with $V(\gf)=E(\gf)=S$, where for the edges
  $s^\alpha=s^\omega=\overline s=s$ hold.
\end{defn}

\begin{lem}
  Let $\gf$ be a graph, and $\gfe=\gf\sqcup\Sigma_\gf$ be the graph
  obtained by adding a loop to every vertex in $\gf$. Let
  $g$ and $e$ be the growth functions for circuits in $\gf$ and $\gfe$
  respectively. Then we have
  $$e(t) = \frac1{1-t}g\left(\frac{t}{1-t}\right).$$
\end{lem}

\begin{defn}[First Product]
  Let $\gfe$ and $\gff$ be two graphs. Their \emph{direct product}
  $\gf=\gfe\times\gff$ is defined by
  $$V(\gf) = V(\gfe)\times V(\gff)$$
  and
  $$E(\gf) = (E(\gfe)\times\Sigma_\gff) \sqcup (\Sigma_\gfe\times E(\gff)).$$
\end{defn}
Note that if the graphs $\gfe$ and $\gff$ have respectively adjacency
matrices $E$ and $F$, then their product has adjacency matrix
$E\tens\idmatrix+\idmatrix\tens F$.

In that case we have
$$G_\gf = \frac1{2i\pi}\oint_{\mathbb S^1}\frac{G_\gfe\left((1+u)t\right)G_\gff\left((1+u^{-1})t\right)}{u} du.$$
This is a simple application of the Laplace transform, that converts
an exponential generating function into an ordinary one and vice
versa~\cite[29.3.3]{abramowitz-s:hbk}.  Indeed, if we had considered
exponential generating functions, the formula would simply have been
$G_\gf = G_\gfe G_\gff$, as is well known~\cite[page 102]{wilf:gf,stanley:gf}.

As an example, let $\gfe=\gff=\Z$, so $G_\gfe = G_\gff = \frac1{\sqrt{1-t^2}}.$
Then
\begin{align*}
  G_\gf &= \frac1{2i\pi}\oint_{\mathbb S^1}\frac{du}{\sqrt{(1-4(1+u)^2t^2)(u^2-4(1+u)^2t^2)}}\\
    &= \frac2\pi K(16t^2) = F\left(\left.\begin{array}{cc}1/2 & 1/2\\ \multicolumn2c{1}\end{array}\right|16t^2\right)
\end{align*}
where $K$ is the complete elliptic function and $F$ the hypergeometric
series.  These functions are known to be transcendental; thus the
circuit series of $\Z^2$ is transcendental. This result appears
in~\cite{grigorchuk-h:groups}.  Numerical evidence suggests the growth function
for $\Z^3$ is not even hypergeometric.

\begin{defn}[Second Product]
  Let $\gfe$ and $\gff$ be two graphs, and suppose that for every
  vertex in $\gfe$ and $\gff$ there is a loop at it. Then their
  \emph{direct product} $\gf=\gfe\times\gff$ is defined by
  $$V(\gf) = V(\gfe)\times V(\gff)$$
  and
  $$E(\gf) = E(\gfe)\times E(\gff).$$
\end{defn}
Note that if the graphs $\gfe$ and $\gff$ have respectively adjacency
matrices $E$ and $F$, then their product has adjacency matrix
$E\tens F$.

In that case we have, again using Laplace transformations
$$G_\gf(t) = \frac1{2i\pi}\oint_{\mathbb S^1}\frac{G_\gfe(u)G_\gff(t/u)}{u} du.$$

Note that with both definitions of products it is possible that the
growth function for circuits in the product be transcendental even if the
growth functions for circuits in the factors are algebraic.

\section{Further Work}
It was mentioned in Subsection~\ref{sec:lang} how the main result applies
to languages. This a special case of a much more general problem:
\begin{problem}
  Given a language $L$ and a set $\mathcal U$ of words, define the
  \emph{desiccation} $L_{\mathcal U}$ of $L$ as the set of words in $L$
  containing no $u\in\mathcal U$ as a subword.
  
  Give sufficient conditions on $L$ and $\mathcal U$ such that a
  formula exist between $\Theta(L)$ and $\Theta\left(L_{\mathcal U}\right)$.
\end{problem}
The special case we studied in this paper is that of
$$\mathcal U=\{s\overline s|\,s\in S\}$$
and a sufficient condition is that $L$ be saturated.

For general $\mathcal U$ this is not always sufficient: let
$S=\{a,b\}$ and $L=b^*(ab^*ab^*)^*$ be the set of words with an even
number of $a$'s. Then if $\mathcal U=\{a^2\}$ there are $7$ desiccated words
of length $5$: $\{b^5, ab^3a, ab^2ab, abab^2, bab^2a, babab,
b^2aba\}$; and if $\mathcal U'=\{b^2\}$ there are $6$ desiccated words of
length $5$: $\{babab, ba^4, aba^3, a^2ba^2, a^3ba, a^4b\}$.
The growth series of $\mathcal U$ and $\mathcal U'$ are the same,
namely $t^2$, but the growth series of $L_{\mathcal U}$ and
$L_{\mathcal U'}$ differ in their degree-$5$ coefficient.

We gave in Section~\ref{sec:freeamalg} a formula relating the circuit
series of a free product to the circuit series of its factors. There
is a notion of \emph{amalgamated product} of graphs, that is a direct
analogue of the amalgamated product of groups.
\begin{problem}
  What conditions on $\gfd,\gfe,\gff$ are sufficient so that
  $$\frac1{(zG_\gf)^{-1}} = \frac1{(zG_\gfe)^{-1}} + \frac1{(zG_\gff)^{-1}} - \frac1{(zG_\gfd)^{-1}}$$
  where $\gf=\gfe*_\gfd\gff$ is an amalgamated product of $\gfe$ and
  $\gff$ along $\gfd$?
\end{problem}

The formula holds if $\gfd$ is the trivial graph; but it cannot hold
in general. If $\gfe=\gff$ is the ``ladder graph'' described in
Section~\ref{ex:ladder}: the set of points $(i,j)$ with $i\in\Z$ and
$j\in\{0,1\}$, with edges connecting all pairs of vertices at
Euclidean distance $1$, and $\gfd$ is $\Z$, embedded as a pole of the
ladder, then the amalgamated product $\gf=\gfe*_\gfd\gff$ is
isomorphic to $\Z^2$. The circuit series of $\gfd$, $\gfe$ and~$\gff$
have been calculated explicitly and are algebraic. The circuit series
of $\gf$ was shown in Section~\ref{sec:direct} to be transcendental;
so there can exist no algebraic definition of $G_\gf$ in terms of
$G_\gfd$, $G_\gfe$ and~$G_\gff$. However, there exists some relations
between these series, as given by~\cite[Theorem~5.5]{voiculescu:rv}.

Given a graph $\gf$, one can construct a graph $\gf^{(k)}$ on the same
vertex set, and with edge set the set of paths of length $\le k$.
Is there some simple relation between the path series of $\gf$ and of
$\gf^{(k)}$? This could be useful for example to obtain asymptotics
about the cogrowth of a group subject to enlargement of generating
set~\cite{champetier:cogrowth}.

The equation~(\ref{eq:free}) corresponds to Voiculescu's
$R$-transform~\cite{voiculescu:rv}. His $S$-transform, in terms of graphs,
corresponds to $\gfe*\gff$ with as edge set all sequences $(e,f)$ and
$(f,e)$, for $e\in E(\gfe)$ and $f\in E(\gff)$. Is there an analogue to
Theorem~\ref{thm:free} in this context?

Finally, (\ref{eq:free}) computes the circuit series of a free product
in terms of the circuit series of the factors. A more complicated
formula yields the path series of a free product in terms of the path
series of the factors. Such considerations give another derivation of
the results in Section~\ref{sec:nonfree}.

\section{Acknowledgements}
The main result of this paper was found in Rome thanks to the
nurturing of Tullio Ceccherini and his family whom I thank. Many
people heard or read preliminary often obscure versions and provided
valuable feedback; I am grateful to (in order of appearance) Michel
Kervaire, Shalom Eliahou, Pierre de la Harpe, Fabrice Liardet,
Rostislav Grigorchuk, Alain Valette, \'Etienne Ghys, Igor Lysionok,
Jean-Paul Allouche, Gilles Robert, and Thierry Vust for their help.

\bibliographystyle{amsalpha}
\providecommand{\bysame}{\leavevmode\hbox to3em{\hrulefill}\thinspace}

\end{document}